\documentclass{article}
\usepackage{latexsym}
\usepackage{amsmath,amssymb,amsthm}

\newtheorem{theorem}{Theorem}[section]
\newtheorem{corollary}[theorem]{Corollary}

\newtheorem{lemma}[theorem]{Lemma}
\newtheorem{remark}[theorem]{Remark}
 \makeatletter
\@addtoreset{equation}{section}

\makeatother

\newcommand{\app}[4]{F_{\!#1}\!
  \left(\left.{#2 \atop #3}\right| #4 \right) }
  \newcommand{\appm}[5]{F^{#1}_{\!#2}\!
  \left(\left.{#3 \atop #4}\right| #5 \right) }

\newcommand{\hpg}[5]{{}_{#1}\mbox{\rm F}_{\!#2}\!
  \left(\left.{#3 \atop #4}\right| #5 \right) }

\newcommand{\hpgo}[2]{{}_{#1}\mbox{\rm F}_{\!#2}}

\newcommand{\equal}{&\!\!\!=\!\!\!&}

\newcommand{\cc}{\lambda}

\newcommand{\RR}{{\Bbb R}}
\newcommand{\CC}{{\Bbb C}}

\newcommand{\PP}{{\Bbb P}}
\newcommand{\ZZ}{{\Bbb Z}}

\title{Dihedral Gauss hypergeometric functions}

\author{Raimundas Vidunas\\
\em Kobe University}

\begin{document}

\maketitle

\begin{abstract} 
Gauss hypergeometric functions with a dihedral monodromy group can be expressed 
as elementary functions, since their hypergeometric equations can be transformed to Fuchsian equations 
with cyclic monodromy groups by a quadratic change of the argument variable.
The paper presents general elementary expressions of these dihedral hypergeometric 
functions, involving finite bivariate sums expressible as terminating Appell's $F_2$ or $F_3$ series.
Additionally, trigonometric expressions for the dihedral functions are presented,
and degenerate cases (logarithmic, or with the monodromy group $\ZZ/2\ZZ$) are considered.
\end{abstract}

\section{Introduction}
\label{logarithms}

As well known, special cases of the Gauss hypergeometric function
$\hpg{2}{1}{A,\,B}{C}{z\,}$ can be represented in terms of elementary functions.
A particularly interesting case are hypergeometric functions with a dihedral monodromy group; 
they can be expressed with square roots inside %compositions of a square root function with
power or logarithmic functions. The simplest examples are:
\begin{eqnarray}   \label{dihedr1}
\hpg{2}{1}{\frac{a}{2},\,\frac{a+1}{2}\,}{a+1}{\,z} \equal
\left(\frac{1+\sqrt{1-z}}{2} \right)^{-a},\\  \label{dihedr2}
\hpg{2}{1}{\frac{a}{2},\,\frac{a+1}{2}\,}{\frac{1}{2}}{\,z} \equal
\frac{(1-\sqrt{z})^{-a}+(1+\sqrt{z})^{-a}}{2},\\ \label{dihedr3}
\hpg{2}{1}{\frac{a+1}{2},\,\frac{a+2}{2}}{\frac{3}{2}}{\,z} \equal
\frac{(1-\sqrt{z})^{-a}-(1+\sqrt{z})^{-a}}{2\,a\,\sqrt{z}} \quad (a\neq 0),\\ 
\label{dihedr4} \hpg{2}{1}{\,\frac{1}{2},\;1\,}{\frac{3}{2}}{\,z} \equal
\frac{\log(1+\sqrt{z})-\log(1-\sqrt{z})}{2\,\sqrt{z}}=\frac{\arctan\sqrt{-z}}{\sqrt{-z}},\\
\label{dihedr5} \hpg{2}{1}{\,\frac{1}{2},\;\frac{1}{2}\,}{\frac{3}{2}}{\,z} \equal
\frac{\log(\sqrt{1-z}+\sqrt{-z})}{\sqrt{-z}}=
\frac{\arcsin\sqrt{z}}{\sqrt{z}}.
\end{eqnarray}
These are solutions of the hypergeometric differential equation 
with the local exponent differences $1/2,1/2,a$ at the three singular points. 
The monodromy group is an infinite dihedral group (for general $a\in\CC$),
or a finite dihedral group (for rational non-integer $a$), 
or an order 2 group (for non-zero integers $a$). 
We refer to Gauss hypergeometric functions with a dihedral monodromy group
as {\em dihedral hypergeometric functions}.

Despite a rich history of research of Gauss hypergeometric functions,
only simplest formulas for dihedral $\hpgo21$ functions like above 
are given in common literature \cite[15.1]{abrostegun}, \cite[2.8]{bateman}. 
General dihedral $\hpgo21$ functions are contiguous to the
simplest functions given above (or their fractional-linear transformations);
that is, their upper and lower parameters differ by integers from the
respective parameters of the simplest functions. The local exponent differences
of their hypergeometric equations are $k+1/2$, $\ell+1/2$ and $\cc\in\CC$, 
where $k,\ell$ are integers. % and $\cc\in\CC$.

An interesting problem is to find explicit elementary expressions for general
dihedral Gauss hypergeometric functions. A set of these expressions is
presented in Sections \ref{sec:explicit}, \ref{sec:easycase}, \ref{sec:trigon} of this paper.
The most general canonical expressions are given in Section \ref{sec:explicit},
in terms of terminating Appell's $F_2$ double sums. 
The key observation is that a particular univariate specialization of Appell's
$F_2$ function satisfies the same Fuchsian equation as a quadratic transformation
of general hypergeometric equations; see Theorem \ref{th:f2gauss}. 
When the monodromy group of the hypergeometric equation is dihedral,
the alluded $F_2$ function is a terminating double sum. Linear relations between
dihedral $\hpgo21$ %Gauss hypergeometric 
and terminating $F_2$ solutions of the same Fuchsian equation give the announced elementary 
expressions for the former. 

Section \ref{sec:easycase} looks at the simpler case of dihedral hypergeometric equations
with the local monodromy differences $k+1/2,1/2,\cc$. % or $k+1/2,1/2,\cc$. 
Then the terminating $F_2$ doublesums become terminating $\hpgo21$ sums. 
The case of the local exponent differences $k+1/2,k+1/2,\cc$ can be reduced 
to the described case by a quadratic transformation of the hypergeometric equation.

Section \ref{sec:degenerate} considers comprehensively
degenerate cases with $\cc\in\ZZ$. 

Section \ref{sec:trigon} considers trigonometric formulas
for dihedral Gauss hypergeometric functions. The simplest such expressions are
\cite[15.1]{abrostegun}, \cite[2.8]{bateman}:
\begin{eqnarray} \label{er:triga}
\hpg21{\,\frac{a}2,\,-\frac{a}2\,}{\frac12}{\sin^2 x} \equal \cos ax,\\ \label{er:trigb}
\hpg21{\frac{1+a}2,\frac{1-a}2}{\frac32}{\sin^2 x} \equal \frac{\sin ax}{a\sin x}.
\end{eqnarray}

Supplementing paper \cite{dihedraltr} presents similar expressions
(as double hypergeometric sums) for quadratic invariants for hypergeometric
equations with the dihedral monodromy group, and describes rational
pull-back transformations between dihedral hypergeometric functions.

\section{Preliminary facts}

Basic facts on hypergeometric functions, Fuchsian equations, the monodromy group,
contiguous relations, terminating hypergeometric sums, Zeilberger's algorithm are well known.
We suggest \cite{specfaar}, \cite{beukers}, \cite{zeilb} as standard though overwhelming references.
Wikipedia pages \cite{wikipedia} can be satisfactorily consulted for quick references.

\subsection{The hypergeometric equation}

The Gauss hypergeometric function $\displaystyle\hpg21{A,\,B}{C}{z}$ 
satisfies the hypergeometric differential equation \cite[Formula (2.3.5)]{specfaar}:
\begin{equation} \label{hpgde}
z\,(1-z)\,\frac{d^2y(z)}{dz^2}+
\big(C-(A+B+1)\,z\big)\,\frac{dy(z)}{dz}-A\,B\,y(z)=0.
\end{equation}
This is a canonical Fuchsian equation on $\PP^1$ %on the complex projective line $\PP^1$
with three singular points. The singularities are $z=0,1,\infty$, % on $\PP^1$.
and the local exponents are:
\[
\mbox{$0$, $1-C$ at $z=0$;} \qquad \mbox{$0$, $C-A-B$ at $z=1$;}\qquad
\mbox{and $A$, $B$ at $z=\infty$.}
\]
The local exponent differences at the singular points are equal (up to a sign) to $1-C$,
\mbox{$C-A-B$} and $A-B$, respectively. Let us denote a hypergeometric equation 
with the local exponent differences $d_1,d_2,d_3$ by $E(d_1,d_2,d_3)$, 
and consider the order of the three arguments unimportant.

Because of frequent use, we recall Euler's and Pfaff's
fractional-linear transformations \cite[Theorem 2.2.5]{specfaar}:
\begin{eqnarray} \label{flinear1}
\hpg{2}{1}{a,\,b\,}{c}{\,z} \equal
(1-z)^{c-a-b}\;\hpg{2}{1}{c-a,\,c-b}{c}{\,z} \\
\label{flinear3} \equal (1-z)^{-a}\;\hpg{2}{1}{a,\,c-b\,}{c}{\frac{z}{z-1}}.
\end{eqnarray}

The following quadratic transformation \cite[(3.1.3),(3.1.9),(3.1.7)]{specfaar}
will illustrate %is used for 
the key reduction of the dihedral monodromy group to a cyclic monodromy group:
 \begin{eqnarray} \label{quadr1}
\hpg{2}{1}{a,\,b}{\frac{a+b+1}{2}}{x} \equal
\hpg{2}{1}{\frac{a}{2},\,\frac{b}{2}}{\frac{a+b+1}{2}}{4x\,(1-x)},\\ \label{quadr2}
\hpg{2}{1}{a,\,b}{a-b+1}{x} \equal (1+x)^{-a}\,
\hpg{2}{1}{\frac{a}{2},\,\frac{a+1}{2}}{a-b+1}{\frac{4x}{(1+x)^2}},\\ \label{quadr3}
\hpg{2}{1}{a,\,b}{2b}{x} \equal \left(1-\frac{x}2\right)^{-a}
\hpg{2}{1}{\frac{a}{2},\,\frac{a+1}{2}}{b+\frac12}{\frac{x^2}{(2-x)^2}}.
\end{eqnarray}
The hypergeometric equations are related in the same way,
by a quadratic {\em pull-back transformation} of the form
\begin{equation} \label{algtransf}
z\longmapsto\varphi(x), \qquad y(z)\longmapsto
Y(x)=\theta(x)\,y(\varphi(x)),
\end{equation}
where $\varphi(x)$ is a rational (quadratic in this case) function, 
and  $\theta(x)$ is a power factor. Geometrically, the transformation {\em pull-backs} 
the starting differential equation on the projective line $\PP^1_z$ to a
differential equation on the projective line $\PP^1_x$, with
respect to the covering $\varphi:\PP^1_x\to\PP^1_z$.
The factor $\theta(x)$ shifts the local exponents of the pull-backed equation,
but it does not change the local exponent differences.
We write the quadratic transformation of hypergeometric equations as
\begin{equation} \label{eq:quadrtr}  \textstyle
E\left( \frac12,\cc,\mu \right) \stackrel2\longleftarrow E\left( \cc,\cc,2\mu \right),
\end{equation}
as the local exponent differences $1/2,\mu$ are doubled by the quadratic substitution,
and the point with the local exponent difference $1/2$ becomes non-singular after
an appropriate choice of $\theta(x)$. The arrow follows the direction of the pull-back covering,
as in \cite{algtgauss}, \cite{dihedraltr}.

\subsection{The monodromy group}

The simplest hypergeometric equations with a dihedral monodromy group are $E(1/2,1/2,a)$.
The monodromy representation for these equations group can be computed using explicit expressions 
(\ref{dihedr1})--(\ref{dihedr4}). If $a\neq 0$, we take (\ref{dihedr2}) and $a\sqrt{z}$ times (\ref{dihedr3})
as a basis of solutions; analytic continuation along loops around $z=0$ and $z=1$ gives
the following generators of the monodromy group:
\begin{equation}
\left( \begin{array}{cc} 1 & 0 \\  0 & -1 \end{array} \right) \qquad\mbox{and}\qquad
\left( \begin{array}{cc} \frac{1+w}2 & \frac{1-w}2 \\  \frac{1-w}2 & \frac{1+w}2 \end{array} \right),
\quad w=\exp(-2\pi i a). 
\end{equation}
If $a=0$, the monodromy generators are 
$\left( \begin{array}{cc} 1 & 0 \\  0 & -1 \end{array} \right)$ and
$\left( \begin{array}{cc} 1 & 2\pi i \\  0 & 1 \end{array} \right)$.

In general, dihedral hypergeometric functions are contiguous to $E(1/2,1/2,a)$.
They are characterized by the property that their differences of local exponents 
at two of the three singular points are half-integers.
If the third local exponent is an irrational number, the monodromy group is an infinite
dihedral group; if it is a non-integer rational number, the monodromy group is a finite
dihedral group. If the third local exponent difference is an integer, the monodromy group
is isomorphic either to $\ZZ/2\ZZ$ %\subset\GG_m$ 
or (in presence of logarithmic solutions) to an infinite dihedral group.
The distinction of these two cases is given in Theorem \ref{th:logc2moodr} below.

The main results of this paper are stated for solutions of 
hypergeometric equation (\ref{hpgde}) with 
\begin{equation} \label{eq:context}
A=\frac{a}{2},\qquad B=\frac{a+1}{2}+\ell,\qquad C=\frac{1}{2}-k.
\end{equation}
The local exponent differences of our working hypergeometric equation are 
$k+1/2$, $\ell+1/2$ and $\cc=a+k+\ell$. Throughout the paper, $k,\ell,m$ denote non-negative 
integers. Except in Section \ref{sec:degenerate}, we assume that $a$ % (or $\cc$) 
is not an integer.

\subsection{The main observation}

If we apply a quadratic pull-back transformation to a hypergeometric equation,
and the two ramified points are singularities of the equation,
the transformed equation is a Fuchsian equation with generally 4 singular points
(and can be solved in terms of {\em Heun functions}). 
Suppose that the hypergeometric equation has a dihedral monodromy group,
hence it is $E\left( k+\frac12,\ell+\frac12,\cc \right)$,
and the two ramified points are the points with the half-integer local exponent differences.
Then the local exponent differences of the transformed equation are $2k+1,2\ell+1,\cc,\cc$.
With an appropriate choice of the power factor $\theta(x)$ in (\ref{algtransf}), 
the transformed equation has trivial monodromy around the two points 
with the integer local exponent differences $2k+1$, $2\ell+1$. 
These two points are not logarithmic because the corresponding points with 
half-integer local exponent differences are not logarithmic.
The monodromy action for the transformed equation will come only from the other two points.
The global monodromy group is therefore cyclic, and the monodromy representation is reducible. 
%the monodromy group conjugate to $\GG_m$ or (in the logarithmic case)  to $\GG_a$. 

In the simplest case $k=0$, $\ell=0$ the transformed equation has just two singularities.
Correspondingly, the classical quadratic transformation (\ref{quadr1}) gives
\begin{equation} \label{fpdihedr}
\hpg{2}{1}{\frac{a}{2},\,\frac{a+1}{2}}{a+1}{4x(1-x)}=
\hpg{2}{1}{a,\,a+1}{a+1}{x}=(1-x)^{-a}.
\end{equation}
This formula is equivalent to (\ref{dihedr1}). If exactly one of $k,\ell$ is zero,
the transformed equation is equivalent to a hypergeometric equation again.
The transformation is then
\begin{equation} \label{eq:quatriv}  \textstyle 
E\left( \frac12,k+\frac12,\cc \right) \stackrel2\longleftarrow E\left( 2k+1,\cc,\cc \right),
\end{equation}
and the transformed solutions can
be expressed in terms of terminating $\hpgo21$ sums, 
as we demonstrate in Section \ref{sec:easycase}.

Even with general integer $k$, $\ell$, the transformed solution must have elementary power 
(at worst logarithmic) solutions, because of the reducible monodromy group.
With a proper normalization by $\theta(x)$ in (\ref{algtransf}) the elementary power solutions
can be polynomials in $x$. It turns out that those polynomials can be written as
terminating Appell's $F_2$ or $F_3$ hypergeometric sums.
Recall that Appell's $F_2$ and $F_3$ bivariate functions are defined by the series
\begin{eqnarray} \label{appf2}
\app2{a;\;b_1,b_2}{c_1,c_2}{x,\,y} \equal \sum_{p=0}^{\infty} \sum_{q=0}^{\infty}
\frac{(a)_{p+q}\,(b_1)_p\,(b_2)_q}{(c_1)_p\,(c_2)_q\;p!\,q!}\,x^p\,y^q,\\ \label{appf3}
\app3{\!a_1,a_2;\,b_1,b_2}{c}{x,\,y} \equal \sum_{i=0}^{\infty} \sum_{j=0}^{\infty}
\frac{(a_1)_p(a_2)_q(b_1)_p(b_2)_q}{(c)_{p+q}\;p!\,q!}\,x^p\,y^q.
\end{eqnarray}
As usual, $(a)_n$ denotes the {\em Pochhammer symbol}  (also called {\em raising factorial}),
which is the product $a(a+1)\cdots(a+n-1)$. 

The key fact is that quadratic transformations of general hypergeometric equations
% (including those with a dihedral monodromy group) 
coincide with the differential equations
for some univariate specializations of Appell's $F_2$ or $F_3$ functions.
In particular, the following is proved in \cite{AppelGhpg}. %Note the 3 free parameters $a,b_1,b_2$.
\begin{theorem} \label{th:f2gauss}
The functions
\begin{equation}\label{eq:f2gauss}
\app2{a;\;b_1,b_2}{2b_1,2b_2}{x,2-x}  \quad\mbox{and}\quad 
(x-2)^{-a}\,\hpg21{\frac{a}2,\frac{a+1}2-b_2}{b_1+\frac12}{\frac{x^2}{(2-x)^2}}
\end{equation}
satisfy the same second order Fuchsian equation.
\end{theorem}
\begin{proof} Part 1 of Theorem 2.4 in \cite{AppelGhpg}. 
\end{proof}
For general values of the 3 parameters $a,b_1,b_2$, the double series for the $F_2(x,2-x)$ function
does not converge for any $x$. However, when $b_1$ and $b_2$ are zero or negative integers,
the $F_2(x,2-x)$ function can be seen as a finite sum of $(1-b_1)(1-b_2)$ terms;
see Remark \ref{rm:simple} below.
On the other hand, the $\hpgo21$ function in (\ref{eq:f2gauss}) is contiguous to 
the $\hpgo21$ function in (\ref{dihedr2}) for integer values of $b_1$, $b_2$, 
hence in general it has a dihedral monodromy group as well. 
This relation between terminating $F_2(x,2-x)$ sums
and dihedral hypergeometric functions is behind our explicit expressions for 
general dihedral functions. We formulate the following variation of Theorem \ref{th:f2gauss}.
\begin{corollary} \label{th:dihf2}
For non-positive integers $k$, $\ell$, the functions
\begin{equation} \label{eq:dih}
\hpg{2}{1}{\frac{a}{2},\,\frac{a+1}{2}+\ell\,}{\frac{1}{2}-k}{\,z\,} \quad\mbox{and}\quad
(1+\sqrt{z})^{-a}\,\app2{a; -k,-\ell}{-2k,-2\ell}{\frac{2\sqrt{z}}{1+\sqrt{z}},\frac2{1+\sqrt{z}}}
\end{equation}
satisfy the same second order Fuchsian equation.
\end{corollary}
\begin{proof}
Substitute $b_1=-k$, $b_2=-\ell$ and $x=2\sqrt{z}/(1+\sqrt{z})$ in Theorem \ref{th:f2gauss}.
\end{proof}
 
We present explicit consequences of this coincidence of differential equations
in Section \ref{sec:explicit}. Note that if $k=0$ or $\ell=0$ then
the double $F_2(x,2-x)$ sum becomes a terminating $\hpgo21(x)$ sum,
in agreement with an observation after formula (\ref{fpdihedr}).

Appell's $F_2(x,y)$ and $F_3(x,y)$ functions are closely related. 
They satisfy the same system of partial differential equations up to a simple transformation,
and particularly, terminating $F_2$ sums become terminating $F_3$ sums when 
summation is reversed in both directions. In particular, for $(a)_{k+\ell}\neq 0$ we have
\begin{eqnarray} \label{eq:r2f3rel}
\app2{a; -k,-\ell}{-2k,-2\ell}{x,y}=\frac{k!\,\ell!\,(a)_{k+\ell}}{(2k)!\,(2\ell)!}\,x^k\,y^{\ell}\,
\app3{k+1,\ell+1; -k,-\ell}{1-a-k-\ell}{\frac1x,\frac1y}.
\end{eqnarray}

\begin{remark} \rm \label{rm:simple}
The hypergeometric series $\displaystyle\hpg{2}{1}{\!-k,\,a}{-2k}{x}$ is not conventionally
defined for a non-negative integer $k$, because of the zero or negative lower parameter.
But it can be usefully interpreted in two ways: as a terminating sum of $k+1$
hypergeometric terms, or by taking the term-wise limit with $k\in\RR$ 
approaching a non-negative positive integer. With both interpretations,
the $\hpgo21$ function is a solution of the respective hypergeometric equation.

In this paper, we adopt the terminating sum interpretation for such $\hpgo21$ functions
and similar bivariate hypergeometric sums. In particular, the $F_2$ function in  (\ref{eq:dih})
is a sum of $(k+1)(\ell+1)$ hypergeometric terms.

% Similarly, the $F_2$ series in (\ref{eq:dih}) and the $F^{2:1;1}_{1:1;1}$ series in (\ref{eq:reverse32}) 
% have the interpretation of terminating summation in both directions (which is the one we use),
% interpretations of terminating summation in one direction but non-terminating summation in the other
% direction, and the interpretation of non-terminating summation in both directions. 
% Either interpretation satisfies the same defining system of partial differential equations.
\end{remark}

\section{Explicit expressions for dihedral functions}
\label{sec:explicit}

The following theorem presents generalizations of (\ref{dihedr1})--(\ref{dihedr3}).
The identities are finite elementary expressions for general dihedral hypergeometric functions. 
The $F_2$ and $F_3$ series are %rational expressions in $\sqrt{z}$, 
finite sums of $(k+1)(\ell+1)$ terms. Because these formulas express solutions
of any dihedral hypergeometric equation, we refer to them as canonical.

Note that the $F_3$ sum in (\ref{eq:diha}) terminates for all (positive or negative) integers $k,\ell$,
as the set of upper parameters does not change under the substitutions $k\mapsto-k-1$ 
and $\ell\mapsto-\ell-1$. 
\begin{theorem}
The following formulas hold for non-negative integers $k,\ell$ and general \mbox{$a\in\CC$}:
\begin{eqnarray} 
\label{eq:diha} \hspace{-9pt}
\hpg{2}{1}{\frac{a}{2},\,\frac{a+1}{2}+\ell}{a+k+\ell+1}{1-z} \!\!\equal\!
z^{k/2} \left(\frac{1+\sqrt{z}}{2} \right)^{-a-k-\ell}
\times\nonumber\\
&&\app3{k+1,\ell+1; -k,-\ell}{a+k+\ell+1}{\frac{\sqrt{z}-1}{2\sqrt{z}},\frac{1-\sqrt{z}}2},\\
\label{eq:dih12} \hspace{-9pt}
\frac{\left(\frac{a+1}2\right)_\ell}{\left(\frac12\right)_\ell}\,
\hpg{2}{1}{\frac{a}{2},\frac{a+1}{2}+\ell}{\frac{1}{2}-k}{\,z} \!\!\equal\! \frac{(1+\sqrt{z})^{-a}}2
\app2{a; -k,-\ell}{-2k,-2\ell}{\frac{2\sqrt{z}}{1+\sqrt{z}},\frac2{1+\sqrt{z}}}\nonumber\\
&& +\frac{(1-\sqrt{z})^{-a}}2
\app2{a; -k,-\ell}{-2k,-2\ell}{\frac{2\sqrt{z}}{\sqrt{z}-1},\frac2{1-\sqrt{z}}},\\ 
\label{eq:dih32} \hspace{-11pt}
\frac{\left(\frac{a+1}2\right)_k\left(\frac{a}2\right)_{k+\ell+1}}
{\left(\frac12\right)_k\left(\frac12\right)_{k+1}\left(\frac12\right)_\ell}\,(-1)^k z^{k+\frac12}\,
\hpg{2}{1}{\frac{a+1}{2}+k,\,\frac{a}{2}+k+\ell+1\,}{\frac{3}{2}+k}{\,z}\hspace{-152pt}\nonumber\\
\!\!\equal\! \frac{(1-\sqrt{z})^{-a}}2
\app2{a; -k,-\ell}{-2k,-2\ell}{\frac{2\sqrt{z}}{\sqrt{z}-1},\frac2{1-\sqrt{z}}}\nonumber\\
&& -\frac{(1+\sqrt{z})^{-a}}2
\app2{a; -k,-\ell}{-2k,-2\ell}{\frac{2\sqrt{z}}{1+\sqrt{z}},\frac2{1+\sqrt{z}}}.
\end{eqnarray}
\end{theorem}
\begin{proof} The first identity follows directly from Karlsson's identity \cite[9.4.(90)]{Srivastava85};
the series on both sides coincide in a neighborhood of $z=1$.
Other local solution at $z=1$ is 
\begin{equation} \label{eq:dihb}
(1-z)^{-a-k-\ell}\,\hpg{2}{1}{-\frac{a}{2}-k-\ell,\,\frac{1-a}{2}-k}{1-a-k-\ell}{1-z},
\end{equation}
which evaluates to
\begin{equation}
z^{k/2} \left(2-2\sqrt{z} \right)^{-a-k-\ell}\,
\app3{k+1,\ell+1; -k,-\ell}{1-a-k-\ell}{\frac{\sqrt{z}-1}{2\sqrt{z}},\frac{1-\sqrt{z}}2}
\end{equation}
by (\ref{eq:diha}) with the substitution $a\mapsto-a-2k-2\ell$. 

The three terms in (\ref{eq:dih12}) are solutions of the same second order
Fuchsian equation by Corollary \ref{th:dihf2}, so there must be  a linear relation between them.
Up to a scalar multiple, the right-hand side of (\ref{eq:dih12}) is the only linear combination 
of the two $F_2$ terms which is invariant under the conjugation $\sqrt{z}\mapsto-\sqrt{z}$. 
Hence the left and right-hand sides of (\ref{eq:dih12}) differ by a factor independent of $z$.
Evaluation of the right side at $z=0$ leads to a terminating $\hpgo21(2)$ sum.
It remains to prove that
\begin{equation} \label{eq:zeilbid}
\hpg21{a,-\ell\,}{-2\ell}{\,2}=\frac{\left(\frac{a+1}2\right)_\ell}{\left(\frac12\right)_\ell}.
\end{equation}
Zeilberger's algorithm is by now a routine technique to
find two-term hypergeometric identities or recurrence relations for hypergeometric sums 
\cite{zeilb}, \cite[Section 3.11]{specfaar}. Its whole output includes {\em certificate} information,
that allows to check the identity or recurrence relation without computer assistance.
Particularly, Zeilberger's algorithm returns the following
difference equation for the $\hpgo21(2)$ summand $S(\ell,j)$:
\begin{equation} \label{eq:zeilb}
(2\ell+1)\,S(\ell+1,j)-(a+1+2\ell)\,S(\ell,j)=H(\ell,j+1)-H(\ell,j),
\end{equation}
where
\begin{equation*}
S(\ell,j)=\frac{(a)_j\,(-\ell)_j}{j!\,(-2\ell)_j}\,2^j, \qquad 
H(\ell,j)=-\frac{j\,(2\ell+1-j)}{2\,(\ell+1-j)}S(\ell,j). % =(1-a-j)S(\ell,j-1).
\end{equation*}
This {\em certificate identity} can be checked by hand straightforwardly. 
To avoid the $0/0$ indeterminancy for $j=\ell+1$, we may also write
$H(\ell,j)=(1-a-j)\,S(\ell,j-1)$. 
While summing up both sides of (\ref{eq:zeilb}) over $j=0,1,\ldots,\ell+1$ 
we obtain a {\em telescoping sum} on the right hand side that simplifies to $H(\ell,\ell+2)-H(\ell,0)=0$.
The summation of the left-hand side gives a first order difference equation
for the $\hpgo21(2)$ sum. The same difference equation is satisfied by 
the right-hand side of (\ref{eq:zeilbid}), and 
a check that the initial values at $\ell=0$ coincide completes the proof of 
(\ref{eq:zeilbid}). % and (\ref{eq:dih12}).

As an intermediate step between formulas (\ref{eq:diha})--(\ref{eq:dih12}) and (\ref{eq:dih32}),
we derive formula (\ref{eq:f2f2rel}) in Lemma \ref{th:f2f3rel} immediately below. 
Then we use connection formula \cite[(2.3.13)]{specfaar}, with both right-side terms
transformed by Euler's transformation (\ref{flinear1}): %\cite[(2.2.7)]{specfaar} 
\begin{eqnarray*}
z^{k+\frac12}\hpg{2}{1}{\!\frac{a+1}{2}+k,\frac{a}{2}+k+\ell+1}{\frac{3}{2}+k}{z} \!=\!
\frac{\Gamma\!\left(\frac32+k\right)\Gamma(-a-k-\ell)}
{\Gamma\!\left(1-\frac{a}2\right)\Gamma\!\left(\frac{1-a}2-\ell\right)}\,
\hpg{2}{1}{\frac{a}{2},\,\frac{a+1}{2}+\ell}{\!a+k+\ell+1}{1\!-\!z}\;\,\\
+ \frac{\Gamma\left(\frac32+k\right)\Gamma(a+k+\ell)}
{\Gamma\!\left(\frac{a+1}2+k\right)\Gamma\!\left(\frac{a}2+k+\ell+1\right)}\,
(1-z)^{-a-k-\ell}\,\hpg{2}{1}{\!-\frac{a}{2}-k-\ell,\,\frac{1-a}{2}-k}{1-a-k-\ell}{1\!-\!z}\!.
\end{eqnarray*}
After applying evaluation (\ref{eq:diha}) to both right-side terms and using relations (\ref{eq:r2f3rel}), (\ref{eq:f2f2rel}), we get (\ref{eq:dih32}).
\end{proof}

\begin{lemma} \label{th:f2f3rel}
We have the following transformation formulas for the terminating $F_2$ and $F_3$ sums
appearing in \mbox{\rm (\ref{eq:diha})--(\ref{eq:dih32})}:  
\begin{eqnarray}  \label{eq:f2f2rel}
&&\hspace{-20pt} \left(1+\sqrt{z}\right)^{k+\ell}
\app2{a;\,-k,-\ell}{-2k,-2\ell}{\frac{2\sqrt{z}}{1+\sqrt{z}},\frac2{1+\sqrt{z}}}=\nonumber\\
 &&  \frac{(-1)^\ell\left(\frac{a+1}2\right)_\ell}{\left(\frac{a+1}2+k\right)_\ell}
\left(1-\sqrt{z}\right)^{k+\ell}
\app2{-a-2k-2\ell;\,-k,-\ell}{-2k,-2\ell}{\frac{2\sqrt{z}}{\sqrt{z}-1},\frac2{1-\sqrt{z}}}.\qquad \\
\label{eq:f3f3rel} && \hspace{-20pt}
\app3{k+1,\ell+1; -k,-\ell}{a+k+\ell+1}{\frac{\sqrt{z}-1}{2\sqrt{z}},\frac{1-\sqrt{z}}2}=\nonumber\\
&& \frac{(a)_{k+\ell}\left(\frac{a+1}2+k\right)_\ell}{(1+a+k+\ell)_{k+\ell}\left(\frac{a+1}2\right)_\ell}
\,\app3{k+1,\ell+1; -k,-\ell}{1-a-k-\ell}{\frac{\sqrt{z}+1}{2\sqrt{z}},\frac{1+\sqrt{z}}2}.
\end{eqnarray}
\end{lemma}
\begin{proof}
Connection formula \cite[(2.3.13)]{specfaar} gives
\begin{eqnarray*}
\hpg{2}{1}{\frac{a}{2},\frac{a+1}{2}+\ell}{\frac{1}{2}-k}{\,z} \equal %=
\frac{\Gamma\!\left(\frac12-k\right)\Gamma(-a-k-\ell)}
{\Gamma\!\left(\frac{1-a}2-k\right)\Gamma\!\left(-\frac{a}2-k-\ell\right)}\,
\hpg{2}{1}{\frac{a}{2},\,\frac{a+1}{2}+\ell}{a+k+\ell+1}{1-z}\;\,\\
&& \hspace{-50pt} + \frac{\Gamma\left(\frac12-k\right)\Gamma(a+k+\ell)}
{\Gamma\!\left(\frac{a}2\right)\Gamma\!\left(\frac{a+1}2+\ell\right)}\,
(1-z)^{-a-k-\ell}\,\hpg{2}{1}{-\frac{a}{2}-k-\ell,\,\frac{1-a}{2}-k}{1-a-k-\ell}{1-z}.
\end{eqnarray*}
We evaluate both terms on the right-hand side using (\ref{eq:diha}), apply (\ref{eq:r2f3rel}) twice
and compare the whole formula with (\ref{eq:dih12}). 
The functions $(1+\sqrt{z})^{-a}$ and $(1-\sqrt{z})^{-a}$ are algebraically independent in general, 
and identification of respective terms to them gives (\ref{eq:f2f2rel}).

Formula (\ref{eq:f3f3rel}) follows from (\ref{eq:f2f2rel}) via (\ref{eq:r2f3rel}).
\end{proof}

Summarising,
the $\hpgo21$ functions (\ref{eq:dih12}) and (\ref{eq:dih32}) form a local basis of solutions at $x=0$;
the functions in (\ref{eq:diha}) and (\ref{eq:dihb}) form a local basis of solutions at $x=1$;
the functions
\begin{equation}
z^{-\frac{a}2}\,\hpg{2}{1}{\frac{a}{2},\frac{a+1}{2}+k}{\frac{1}{2}-\ell}{\,\frac1z},\qquad
z^{-\frac{a+1}2-\ell}\,\hpg{2}{1}{\frac{a+1}{2}+\ell,\,\frac{a}{2}+k+\ell+1}{\frac{3}{2}+\ell}{\,\frac1z}
\end{equation}
form a local basis of solutions at $x=\infty$. Each of the six functions can be (generally) 
expressed  by a hypergeometric series in 4 different ways following Euler-Pfaff transformations 
(\ref{flinear1})--(\ref{flinear3}). This gives the standard set (and structure) of
$6\times4=24$ Kummer's hypergeometric solutions of hypergeometric equation (\ref{hpgde}). 
The symmetries of  (\ref{hpgde}) correspond to the permutations of the 3 singular points and 
of local exponents at them. Within our context described by (\ref{eq:context}), 
the action of these permutations on the parameters of hypergeometric functions is the following:
\begin{itemize}
\item Permutation of the local exponents at $z=0$: the parameters are transformed as
\mbox{$k\mapsto-k-1$}, $a\mapsto a+2k+1$; the hypergeometric solution gets multiplied by $z^{k+1/2}$.
\item Permutation of the local exponents at $z=1$: the parameters are transformed as
\mbox{$a\mapsto -a-2k-2\ell$}; the hypergeometric solution gets multiplied by $(1-z)^{-a-k-\ell}$.
\item Permutation of the local exponents at $z=\infty$: the parameters are transformed as
\mbox{$\ell\mapsto-\ell-1$}, $a\mapsto a+2\ell+1$; the hypergeometric solution gets multiplied 
by $z^{-\ell-1/2}$.
\item Permutation $z\mapsto 1/z$ of the singularities $z=0$, $z=\infty$: the parameters are 
transformed as $k\leftrightarrow\ell$; the hypergeometric solution gets multiplied by $z^{-a/2}$.
\item Permutation $z\mapsto 1-z$ of the singularities $z=0$, $z=1$: the parameters are 
transformed as $k\mapsto-a-k-\ell-\frac12$.
\item Permutation $z\mapsto z/(z-1)$ of the singularities $z=1$, $z=\infty$: the parameters are 
transformed as $\ell\mapsto-a-k-\ell-\frac12$;  the hypergeometric solution gets multiplied 
by \mbox{$(1-z)^{-a/2}$}.
\end{itemize}
Pfaff's fractional-linear transformation (\ref{flinear3}) %\cite[Theorem 2.2.5]{specfaar} 
gives a different relation between the upper parameters in (\ref{eq:dih12}) or (\ref{eq:dih32});
for example
\begin{eqnarray}
\hpg{2}{1}{\frac{a}{2},\frac{a+1}{2}+\ell}{\frac{1}{2}-k}{\,z}
%\equal (1-z)^{-a-k-\ell}\,\hpg{2}{1}{-\frac{a}{2}-k-\ell,\,\frac{1-a}{2}-k}{\frac{1}{2}-k}{\,z}\\ 
\label{eq:eupf12}
\equal (1-z)^{-\frac{a}2}\,\hpg{2}{1}{\frac{a}{2},-\frac{a}{2}-k-\ell}{\frac{1}{2}-k}{\frac{z}{z-1}}.
%\\ \equal (1-z)^{-\frac{a+1}2-\ell}\,\hpg{2}{1}{\frac{1-a}{2}-k,\frac{a+1}{2}+\ell}{\frac{1}{2}-k}{\frac{z}{z-1}}.
\end{eqnarray}
The alternative shape of the upper parameters is convenient in Section \ref{sec:trigon}.

\section{The simple cases}
\label{sec:easycase}

In the special cases $k=0$ or $\ell=0$, one of the local exponents of dihedral hypergeometric
equation (\ref{hpgde}) is equal to $1/2$. %(that is,  if $k=0$ or $\ell=0$),
Then the terminating double $F_2$ or $F_3$ sums in (\ref{eq:diha})--(\ref{eq:dih32}) become 
terminating single $\hpgo21$ sums. 
For example, special cases of (\ref{eq:dih12})--(\ref{eq:dih32}) are
\begin{eqnarray} \label{eq:diha12a}
&& \hpg{2}{1}{\!\frac{a}{2},\frac{a+1}{2}}{\frac{1}{2}-k}{z} = 
\nonumber  \\ && \hspace{40pt}
\frac{(1\!+\!\sqrt{z})^{-a}}2\hpg{2}{1}{\!-k,a}{-2k}{\frac{2\sqrt{z}}{1\!+\!\sqrt{z}}}  %\nonumber \\ && 
+\frac{(1\!-\!\sqrt{z})^{-a}}2\hpg{2}{1}{\!-k,a}{-2k}{\frac{2\sqrt{z}}{\sqrt{z}\!-\!1}}, \qquad\\
&& \frac{\left(\frac{a+1}2\right)_\ell}{\left(\frac12\right)_\ell}\,
\hpg{2}{1}{\!\frac{a}{2},\frac{a+1}{2}+\ell}{\frac{1}{2}}{z} = 
\nonumber \\ && \hspace{40pt}
\frac{(1\!+\!\sqrt{z})^{-a}}2\hpg{2}{1}{\!-\ell,a}{-2\ell}{\frac{2}{1\!+\!\sqrt{z}}}  %\nonumber \\ && 
+\frac{(1\!-\!\sqrt{z})^{-a}}2\hpg{2}{1}{\!-\ell,a}{-2\ell}{\frac{2}{1-\sqrt{z}}},\\
&& \frac{(-1)^k\,\left(a\right)_{2k+1}}
{2^{2k+1}\left(\frac12\right)_k\left(\frac12\right)_{k+1}}\, z^{k+\frac12}\,
\hpg{2}{1}{\frac{a+1}{2}+k,\,\frac{a}{2}+k+1\,}{\frac{3}{2}+k}{\,z}=
\nonumber  \\ && \hspace{40pt}
\frac{(1\!-\!\sqrt{z})^{-a}}2\hpg{2}{1}{\!-k,a}{-2k}{\frac{2\sqrt{z}}{\sqrt{z}\!-\!1}}  %\nonumber \\ && 
-\frac{(1\!+\!\sqrt{z})^{-a}}2\hpg{2}{1}{\!-k,a}{-2k}{\frac{2\sqrt{z}}{1\!+\!\sqrt{z}}}, \qquad\\
&& \frac{2\left(\frac{a}2\right)_{\ell+1}}{\left(\frac12\right)_\ell}\;\sqrt{z}\,
\hpg{2}{1}{\frac{a+1}{2},\,\frac{a}{2}+\ell+1\,}{\frac{3}{2}}{\,z} =
\nonumber \\ && \hspace{40pt}
\frac{(1\!-\!\sqrt{z})^{-a}}2\hpg{2}{1}{\!-\ell,a}{-2\ell}{\frac{2}{1\!-\!\sqrt{z}}}  %\nonumber \\ && 
-\frac{(1\!+\!\sqrt{z})^{-a}}2\hpg{2}{1}{\!-\ell,a}{-2\ell}{\frac{2}{1\!+\!\sqrt{z}}}.
\end{eqnarray}
This appearance of terminating $\hpgo21$ sums is expectable,
since quadratic transformation (\ref{eq:quatriv})
% (\ref{eq:qled}) becomes \begin{equation} \label  \textstyle 
% E\left( \frac12,k+\frac12,\cc \right) \stackrel2\longleftarrow E\left( 2k+1,\cc,\cc \right),
% \end{equation}
leads to a hypergeometric equation with a cyclic monodromy group. 
We have quadratic transformations between the dihedral and terminating $\hpgo21$ functions.
In particular,  the identification $z=4x/(1+x)^2$ in classical formula (\ref{quadr2}) gives
\begin{equation} \label{eq:fpdihedr}
\hpg{2}{1}{\frac{a}{2},\,\frac{a+1}{2}}{a+k+1}{\,z} 
=\left(\frac{1+\sqrt{1-z}}{2} \right)^{-a}\,\hpg21{-k,\,a}{a+k+1}{\frac{1-\sqrt{1-z}}{1+\sqrt{1-z}}},
\end{equation}
while formulas (\ref{quadr1}) and (\ref{flinear3}) give
\begin{eqnarray} 
\label{eq:fpdihedr2a} \hpg{2}{1}{\frac{a}{2},\,\frac{a+1}{2}+\ell}{a+\ell+1}{\,z} 
\equal \left(\frac{1+\sqrt{1-z}}{2} \right)^{-a}\,\hpg21{-\ell,\,a}{a+\ell+1}{\frac{\sqrt{1-z}-1}{\sqrt{1-z}+1}}.
\end{eqnarray}
We recognize here special cases of (\ref{eq:diha}) after the substitution $z\mapsto 1-z$ 
and application of Pfaff's fractional-linear transformation (\ref{flinear3}) %\cite[(2.2.6)]{specfaar} 
to the terminating $\hpgo21$ series.
Relatedly, formulas (\ref{eq:f2f2rel}), (\ref{eq:f3f3rel}) are standard hypergeometric identities
when $k=0$ or $\ell=0$.

The dihedral functions with $k=\ell$ can be reduced to the considered case $k\,\ell=0$
via a quadratic transformation. The quadratic transformation is
\begin{eqnarray} \label{eq:quadih}  \textstyle
E\left(\frac12,\,k+\frac12,\,\cc\right) \stackrel2\longleftarrow E\left(k+\frac12,\,k+\frac12,\,2\cc \right).
\end{eqnarray}
Quadratic transformations (\ref{quadr2}) and (\ref{quadr3}) give the following expressions:
\begin{eqnarray} \label{eq:symmkl}
&&\hspace{-20pt} \hpg21{a,a+k+\frac12}{\frac12-k}{z} =
(1+z)^{-a}\,\hpg21{\frac{a}2,\,\frac{a+1}2}{\frac12-k}{\frac{4z}{(1+z)^2}}  \nonumber\\
&& = \frac{(1\!+\!\sqrt{z})^{-2a}}2\hpg{2}{1}{\!-k,a}{-2k}{\frac{4\sqrt{z}}{(1\!+\!\sqrt{z})^2}} 
+\frac{(1\!-\!\sqrt{z})^{-2a}}2\hpg{2}{1}{\!-k,a}{-2k}{-\frac{4\sqrt{z}}{(\sqrt{z}\!-\!1)^2}}\!, \qquad\\
&&\hspace{-20pt} \hpg21{a,\,a+k+\frac12}{2a+2k+1}{1-z} =
\left(\frac{1+z}2\right)^{-a}\,\hpg21{\frac{a}2,\frac{a+1}2}{a+k+1}{\frac{(1-z)^2}{(1+z)^2}}
\nonumber \\ &&
=z^{k/2}\left(\frac{1+\sqrt{z}}2\right)^{-2a-2k}\hpg21{-k,k+1}{a+k+1}{-\frac{(\sqrt{z}\!-\!1)^2}{4\sqrt{z}}}.
\end{eqnarray}
Comparing with formulas (\ref{eq:dih12}), (\ref{eq:diha}) for the left-hand sides here, we conclude
\begin{eqnarray}
\app2{2a;\,-k,-k}{-2k,-2k}{x,2-x}\equal 
\frac{\left(a+\frac12\right)_k}{\left(\frac12\right)_k}\,\hpg21{-k,a}{-2k}{x(2-x)}\!,\\
\app3{-k,-k;\,k+1,k+1}{2a+2k+1}{x,\frac{x}{2x-1}}\equal\hpg21{-k,k+1}{a+k+1}{\frac{x^2}{2x-1}}.
\end{eqnarray}
The finite $\hpgo21$ sums can be modified using these transformation formulas 
\cite[Section 7]{degeneratehpg}: 
\begin{eqnarray} \label{eq:termlz}
\hpg{2}{1}{-k,\,a\,}{-2k}{\,x} &=& \left(1-x\right)^k\,\hpg{2}{1}{-k,\,-a-2k}{-2k}{\frac{x}{x-1}}\\
\label{eq:termlzb}
&=& \frac{k!\,(a)_k}{(2k)!}\,x^{k}\,\hpg{2}{1}{-k,\,k+1}{1-a-k}{\,\frac1x} \\
\label{eq:termlzc}
&=& \frac{k!\,(1+a+k)_k}{(2k)!}\,x^k\,\hpg{2}{1}{-k,\,k+1}{1+a+k}{1-\frac1x} \\
&=& \frac{k!\,(1+a+k)_k}{(2k)!}\,\hpg21{-k,\,a}{1+a+k}{1-x} \\ \label{eq:termlzz}
&=& \frac{k!\,(a)_k}{(2k)!}\left(x-1\right)^k\,\hpg21{-k,\,-a-2k}{1-a-k}{\frac{1}{1-x}}.
\end{eqnarray}
These six hypergeometric expressions %in (\ref{eq:termlz})--(\ref{eq:termlzz})
have the following arguments after the substitution \mbox{$x=2\sqrt{z}\big/\left(1+\sqrt{z}\right)$}, respectively:
\begin{eqnarray}
\frac{2\sqrt{z}}{1+\sqrt{z}}, \quad \frac{2\sqrt{z}}{\sqrt{z}-1},\quad
\frac{1+\sqrt{z}}{2\sqrt{z}}, \quad \frac{\sqrt{z}-1}{2\sqrt{z}}, \quad
\frac{\sqrt{z}-1}{\sqrt{z}+1}, \quad \frac{\sqrt{z}+1}{\sqrt{z}-1}.
\end{eqnarray}
It is easy to substitute further $z\mapsto 1/z$ or $z\mapsto 1-z$ here.
The substitution $z\mapsto z/(z-1)$ leads to the following six arguments, respectively:
\begin{eqnarray}
2z-2\sqrt{z^2-z}, \quad 2z+2\sqrt{z^2-z},\quad
\frac12+\frac{\sqrt{z^2-z}}{2z}, \quad \frac12-\frac{\sqrt{z^2-z}}{2z}, \nonumber \\
1-2z+2\sqrt{z^2-z}, \quad 1-2z-2\sqrt{z^2-z}.
\end{eqnarray}
The argument on the right-hand side of (\ref{eq:fpdihedr2a}) 
can be written as
\begin{equation}
\frac{\sqrt{1-z}-1}{\sqrt{1-z}+1}=1-\frac{2}z+\frac{2\sqrt{1-z}}{z}.
\end{equation}

\begin{remark} \rm \label{kummer24}
As well known  \cite[Section 2.9]{specfaar}, there are generally 24 hypergeometric series that are
solutions of the same hypergeometric equation (\ref{hpgde}); they are referred to as
24 {\em Kummer's solutions}. Particularly,  a general hypergeometric equation has
a basis of hypergeometric solutions at each of the three singular points; 
the six solutions are different functions, and each of them has 4 representations
as Gauss hypergeometric series due to Euler-Pfaff transformations (\ref{flinear1})--(\ref{flinear3}).
When terminating or logarithmic solutions are present, this structure of 24 solutions degenerates
\cite{degeneratehpg}.

The hypergeometric equation on the left-hand side of transformation (\ref{eq:quatriv})
has a degenerate structure of 24 Kummer's solutions,
as exemplified by formulas (\ref{eq:termlz})--(\ref{eq:termlzz}).
According to \cite[Section 7]{degeneratehpg},
the same hypergeometric equation has other terminating solution
\begin{equation}  \label{eq:termlz2}
(1-x)^{-a-k}\,\hpg{2}{1}{-k,\,-a-2k}{-2k}{\,x}=(1-x)^{-a}\,\hpg{2}{1}{-k,\;a}{-2k}{\frac{x}{x-1}}.
\end{equation}
When written in terms of $z$ under the identification \mbox{$x=2\sqrt{z}\big/\left(1+\sqrt{z}\right)$},
both terminating solutions are related (up to a power factor) by the conjugation $\sqrt{z}\mapsto-\sqrt{z}$.
The two different terminating solutions are present on the left-hand side of (\ref{eq:diha12a}). 
Both terminating solutions are representable by 6 terminating and 4 non-terminating $\hpgo21$ sums.
The remaining 4 Kummer's solutions of the transformed equation
are non-terminating $\hpgo21$ series at $x=0$. They are related by
Euler-Pfaff transformations (\ref{flinear1})--(\ref{flinear3}), and represent the following solution:
\begin{eqnarray} \label{eq:snonterm}
\frac{(-1)^{k}\,k!^2\,(a)_{2k+1}}{(2k)!\,(2k+1)!}x^{2k+1}\hpg21{k+1,\,a+2k+1}{2k+2}{\,x}=
\nonumber\hspace{100pt}\\
(1-x)^{-a-k}\,\hpg{2}{1}{-k,\,-a-2k}{-2k}{\,x}-\hpg{2}{1}{-k,\,a}{-2k}{\,x}.
\end{eqnarray}
Quadratic transformation (\ref{quadr3}) gives the following identification:
\begin{equation} \hspace{4pt}
\hpg{2}{1}{k+1,a+2k+1}{2k+2}{\frac{2\sqrt{z}}{1+\sqrt{z}}} \!=
\left(1+\sqrt{z}\right)^{a}\hpg{2}{1}{\!\frac{a+1}{2}+k,\,\frac{a}{2}+k+\ell+1}{\frac{3}{2}+k}{z}\!,
\end{equation}
consistent with (\ref{eq:dih32}). Formulas (\ref{flinear1})--(\ref{quadr3}) for fractional-linear
and quadratic transformations may not hold when a degenerate set of 24 Kummer's solutions
is involved. In particular, the terminating $\hpgo21$ sums in (\ref{eq:termlz}) and (\ref{eq:termlz2}) 
are not related by Euler-Pfaff transformations (\ref{flinear1})--(\ref{flinear3}). 
As noted in \cite[Lemma 3.1]{degeneratehpg}, the generic transformations 
are correct here only if exactly one of the $\hpgo21$ sums is interpreted 
as non-terminating (following Remark \ref{rm:simple}). The right-hand side of (\ref{eq:snonterm}) 
can be seen as the difference between
the terminating and non-terminating interpretations of $\displaystyle\hpg{2}{1}{-k,\,a}{-2k}{\,x}$.
Relatedly, there are no two-term identities between the dihedral function in  (\ref{eq:fpdihedr})
and the terminating $\hpgo21$ sums in (\ref{eq:termlz}) or (\ref{eq:termlz2}).

Transformations  (\ref{eq:quatriv}) and (\ref{eq:quadih}) involve degenerate 
Gauss hypergeometric functions on both sides when $\cc$ is an integer.
We comment further on degeneracies of Kummer's 24 solutions 
in Subsection {sc:degenerate} below. % Remarks \ref{rm:degenerate} and \ref{rm:degtrans} 
\end{remark}

\begin{remark} \label{rm:legendre} \rm
The dihedral functions with $k=0$ or $\ell=0$ can be expressed in terms of 
the {\em (associated) Legendre functions} $P_{\nu}^{\mu}(z)$,  $Q_{\nu}^{\mu}(z)$ 
with integer $\nu$ or half-integer $\mu$.
The relation to the Legendre functions with integer $\nu$ is clear after comparing the
terminating series in (\ref{eq:termlzb}) with the definition \cite[8.1.2]{abrostegun}:
\begin{equation}
P_k^{\mu}(x)=\frac1{\Gamma(1-\mu)}\left(\frac{x+1}{x-1}\right)^{\frac{\mu}2}
\hpg21{-\nu,\,\nu+1}{1-\mu}{\frac{1-x}2}.
\end{equation}
By formulas in \cite[Chapter 8]{abrostegun} or \cite[Chapter III]{bateman} we obtain these expressions:
\begin{eqnarray} \label{eq:diha12b}
\frac{(2k-1)!!}{\Gamma(1-a)}\,z^{-\frac{k}2}\,(1-z)^{\frac{a+k}2}\,
\hpg{2}{1}{\!\frac{a}{2},\frac{a+1}{2}}{\frac{1}{2}-k}{z} \!\!\!\equal\!\! 
\frac12\,P_k^{a+k}\!\left(\frac1{\sqrt{z}}\right)
+\frac{(-1)^k}2\,P_k^{a+k}\!\left(-\frac1{\sqrt{z}}\right) \qquad\quad \\
\equal\!\! P_k^{a+k}\!\left(\frac1{\sqrt{z}}\right)
-\frac{\sin\pi a}{(-1)^a\,\pi}\,Q_k^{a+k}\!\left(\frac1{\sqrt{z}}\right),\\
z^{\frac{k+1}2}\,(1-z)^{\frac{a+k}2}\,
\hpg{2}{1}{\!\frac{a+1}{2},\frac{a}{2}+k+1}{\frac{3}{2}+k}{z} \!\!\!\equal\!\! 
\frac{(2k+1)!!}{\Gamma(a+2k+1)}\,(-1)^{-a-k}\,Q_k^{a+k}\!\left(\frac1{\sqrt{z}}\right),\\
\frac{2^{\ell+1}\!\left(\frac{a+1}2\right)_{\ell}(1-z)^{\frac{a+\ell}2}}{(-1)^{\frac{a-\ell}2}\,\Gamma(1-a)}\,
\hpg{2}{1}{\!\frac{a}{2},\frac{a+1}{2}+\ell}{\frac{1}{2}}{\,z} \!\!\!\equal\!\!
P_k^{a+\ell}\!\left(\sqrt{z}\right)+(-1)^{-a-\ell}\,P_k^{a+\ell}\!\left(-\sqrt{z}\right),\\
\frac{2^{\ell+2}\!\left(\frac{a}2\right)_{\ell+1}(1-z)^{\frac{a+\ell}2} \sqrt{z}}
{(-1)^{\frac{a-\ell}2}\,\Gamma(1-a)}\,
\hpg{2}{1}{\!\frac{a+1}{2},\frac{a}{2}+\ell+1}{\frac{3}{2}}{\,z} \hspace{-35pt} \nonumber\\
\equal\!\! P_k^{a+\ell}\!\left(\sqrt{z}\right)-(-1)^{-a-\ell}\,P_k^{a+\ell}\!\left(-\sqrt{z}\right),\\
z^{-\frac{k}2} \left(\frac{1-z}4\right)^{\!\frac{a+k}2}
\hpg21{\frac{a}2,\frac{a+1}2}{a+k+1}{1-z}   \!\!\!\equal\!
\Gamma(1+a+k)\,P_k^{-a-k}\!\left(\frac1{\sqrt{z}}\right).
\end{eqnarray}
We generally mean $(-1)^x=\exp(i\pi x)$ here. For integer $k$, 
formulas \cite[8.2.3, 8.2.5]{abrostegun} give the simpler relation
\begin{equation}
P_k^{-\mu}(x)=(-1)^k\,\frac{\Gamma(1-\mu+k)}{\Gamma(1+\mu+k)}\,P_k^{\mu}(-x),
\end{equation}
while \cite[8.2.4, 8.2.6]{abrostegun} give
\begin{equation}
Q_k^{\mu}(-x)=(-1)^{k+1} Q_k^{\mu}(x), \qquad
Q_k^{-\mu}(x)=(-1)^{-2\mu}\,\frac{\Gamma(1-\mu+k)}{\Gamma(1+\mu+k)}\,Q_k^{\mu}(x).
\end{equation}
Expressions in terms of Legendre functions with half-integer $\mu$ are obtained by applying  
\cite[8.2.7--8]{abrostegun}:
\begin{eqnarray*}
P_{-\mu-\frac12}^{-k-\frac12}\left(\frac1{\sqrt{1-z}}\right) \equal
\frac{(-1)^{-\mu}\,\sqrt{2/\pi}}{\Gamma(\mu+k+1)} \left(\frac{1-z}{z}\right)^{\!\frac14}
Q_k^{\mu}\!\left(\frac1{\sqrt{z}}\right),\\
Q_{-\mu-\frac12}^{-k-\frac12}\left(\frac1{\sqrt{1-z}}\right) \equal
i\,(-1)^{k+1}\Gamma(-\mu-k)\sqrt{\frac{\pi}2} \left(\frac{1-z}{z}\right)^{\!\frac14}
P_k^{\mu}\!\left(\frac1{\sqrt{z}}\right).
\end{eqnarray*}
\end{remark}

\section{Degenerate and logarithmic solutions}
\label{sec:degenerate}

%Let $\GG_m$ and $\GG_a$ denote the representations 
%$\left\{ \left( \begin{array}{cc} 1 & 0 \\  0 & v \end{array} \right)\;:\;v\in\CC^*\right\}$
%and $\left\{ \left( \begin{array}{cc} 1 & v \\  0 & 1 \end{array} \right)\;:\;v\in\CC\right\}$
%of, respectively, the multiplicative group $\CC^*$ and the additive group $\CC$.

So far we considered solutions of hypergeometric equations $E(k+1/2,\ell+1/2,\cc)$, %$a+k+\ell$, 
where $k,\ell$ are integers but $\cc$ is not an integer.
If the third local exponent difference $\cc$ is an integer, the monodromy group
is either completely reducible and isomorphic to $\ZZ/2\ZZ$ %\subset\GG_m$ 
or (in presence of logarithmic solutions) it is isomorphic to an infinite dihedral group.
First we state conditions how to separate the two cases.

\subsection{Conditions for logarithmic solutions}

\begin{theorem} \label{th:logc2moodr}
Let $k,\ell,n$ denote non-negative integers.
Then equation $E(\mbox{$k+1/2$},\mbox{$\ell+1/2$},n)$ has logarithmic solutions
if and only if
\begin{eqnarray} \label{logcond1}
n \le k+\ell,\;\, && \mbox{if $n+k+\ell$ is even},  \\ \label{logcond2}  
n< \left| k-\ell \right|, && \mbox{if $n+k+\ell$ is odd}.   
\end{eqnarray}
If this is the case, the monodromy group of $(\ref{hpgde})$ is an infinite dihedral group;
otherwise the monodromy group is isomorphic to $\ZZ/2\ZZ$. % \subset\GG_m$.
\end{theorem}
\begin{proof} A representative equation (\ref{hpgde}) with the assumed local exponents has 
\[
A=-\frac{n+k+\ell}2, \qquad B=-\frac{n+k-\ell-1}2, \qquad C=\frac12-k.
\] 
The sequence $A$, $1-B$, $C-A$, $1+B-C$ contains exactly two integers.
By part (3) of \cite[Theorem 2.2]{degeneratehpg}, there are no logarithmic solutions
precisely when the two integers are either both positive or both non-positive.
Equivalently, there are logarithmic solutions precisely when one of the integers is
positive while the other is non-positive. If $n+k+\ell$ is even, the two integers are
\[
A=-\frac{n+k+\ell}2 \qquad\mbox{and}\qquad 1+B-C=1+\frac{k+\ell-n}2.
\]
The first integer is always zero or negative; the second integer is positive exactly when 
$n\le k+\ell$. If $n+k+\ell$ is odd, the two integers are
\[
1-B=\frac{1+n+k-\ell}2 \qquad\mbox{and}\qquad C-A=\frac{1+n-k+\ell}2.
\]
We may assume $\ell\le k$ without loss of generality. Then the first integer is positive;
the second integer is non-positive exactly when $n<k-\ell$.
\end{proof}

For comparison, recall \cite[Section 9]{degeneratehpg} 
that hypergeometric equation $E(k,\ell,n)$ with non-negative integers $k,\ell,n$ 
has logarithmic solutions if and only if 
one of the integers is greater than the sum of the other two.
Here is a more direct formulation of Theorem \ref{th:logc2moodr}.
\begin{corollary} %\label{lm:logc2moodr}
Suppose that $p,q$ are half-integers, %$\lb=k+1/2,\lc=\ell+1/2$
and $n$ is a non-negative integer. 
The set \mbox{$\{|p-q|,|p+q|\}$} contains two integers of different parity;
let $K$ be the integer in this set such that $K+n$ is odd.
Then equation $E(p,q,n)$ has the monodromy group isomorphic to $\ZZ/2\ZZ$ if $K<n$,
and it has logarithmic solutions otherwise.
\end{corollary}

Here is a formulation of Theorem \ref{th:logc2moodr} that refers to the parameter $a=\cc+k+\ell$ 
in (\ref{eq:context}) rather than to the local exponent difference $\cc$.
\begin{corollary} \label{lm:logc2moodr}
Let $k,\ell,m$ denote non-negative integers, and suppose that $k\le\ell$. 
Hypergeometric equation $(\ref{hpgde})$ with $(\ref{eq:context})$ and
%for the $\hpgo21$ functions in $(\ref{eq:diha})$--$(\ref{eq:dihb})$ with 
$a=-m$  has logarithmic solutions if and only if
\begin{eqnarray} \label{eq:logcond1}
0\le \frac{m}2 \le k+\ell, && \mbox{for even $m$},  \\ \label{eq:logcond2}  
k<\frac{m+1}2\le \ell, && \mbox{for odd $m$}.   
\end{eqnarray}
If this is the case, the monodromy group of $(\ref{hpgde})$ is an infinite dihedral group;
otherwise the monodromy group is isomorphic to $\ZZ/2\ZZ$. %\subset\GG_m$.
\end{corollary}
\begin{proof} Theorem \ref{th:logc2moodr} is being applied to equation 
$E(k+1/2,\ell+1/2,m-k-\ell)$ or to \mbox{$E(k+1/2,\ell+1/2,k+\ell-m)$}.
\end{proof}

We refer to solutions of a hypergeometric equation with the monodromy group $\ZZ/2\ZZ$ as 
{\em degenerate}. Degenerate or logarithmic solutions are involved exactly when a Pochhammer
factor on the left-hand side of our main identities (\ref{eq:dih12})--(\ref{eq:dih32}) vanishes. 
The identities are still valid even if the left-hand side vanishes.
Recalling $\cc=a-k-\ell$, we set 
\begin{equation}
K_1=\min(k,\ell), \qquad K_2=\ell-k, \qquad L=k+\ell.
\end{equation}
Then Pochhammer factors in both (\ref{eq:dih12}) and (\ref{eq:dih32}) vanish if and only if
\begin{eqnarray} 
\label{eq:conda1} && a\in\{-1,-3,\ldots,1-2K_1\}, \\
\label{eq:conde1} && \cc\in\{|K_2|+1,|K_2|+3,\ldots,L-1\}.
\end{eqnarray}
Only a Pochhammer factor in (\ref{eq:dih32}) vanishes if and only if
\begin{eqnarray}
\label{eq:conda2} && a\in\{0,-2,\ldots,-2L)\}\cup\{-1-2K_1,-3-2K_1,\ldots,1-2k\}, \\
\label{eq:conde2} && \cc\in\{-L,2-L,\ldots,L\}\cup\{K_2+1,K_2+3,\ldots,|K_2|-1\},
\end{eqnarray}
and only the Pochhammer factor in (\ref{eq:dih12}) vanishes if and only if
\begin{eqnarray}
\label{eq:conda3} && a\in\{-1-2K_1,-3-2K_1,\ldots,1-2\ell\},\\
\label{eq:conde3} && \cc\in\{-K_2+1,-K_2+3,\ldots,|K_2|-1\}.
\end{eqnarray}
By Theorem \ref{th:logc2moodr}, hypergeometric equation $E(k+1/2,\mbox{$\ell+1/2$},\cc)$ has logarithmic solutions if and only if exactly one of conditions (\ref{eq:conde2}) or (\ref{eq:conde3}) holds. Equivalently, by Corollary \ref{lm:logc2moodr} we have to check whether exactly one of conditions (\ref{eq:conda2}) or (\ref{eq:conda3}) holds.

\subsection{The monodromy group $\ZZ/2\ZZ$}

By Corollary \ref{lm:logc2moodr}, the monodromy group of hypergeometric equation
$(\ref{hpgde})$ with $(\ref{eq:context})$ and $a=-m$ is isomorphic to $\ZZ/2\ZZ$
if and only if either (\ref{eq:conda1}) holds or both (\ref{eq:conda2}) and (\ref{eq:conda3}) hold.
%If (\ref{eq:conda1}) holds, %according to Corollary \ref{lm:logc2moodr}. 
In the former case,
formulas (\ref{eq:dih12}) and (\ref{eq:dih32}) do not give terminating expressions for the functions
\begin{equation} \label{eq:c2hpg}
\hpg{2}{1}{-\frac{m}{2},\,\ell-\frac{m-1}{2}}{\frac{1}{2}-k}{\,z}, \qquad
\hpg{2}{1}{k-\frac{m-1}{2},\,k+\ell-\frac{m}{2}+1\,}{\frac{3}{2}+k}{z}.
\end{equation}
But terminating expressions for these functions are obtained by applying Euler's transformation
(\ref{flinear1}). Adding up the right-hand sides of  (\ref{eq:dih12}) and (\ref{eq:dih32}) 
suggests the following lemma.
\begin{lemma} \label{th:zerof2}
Suppose that $m$ is an odd positive integer.
Then for any $b,c\in\CC$ 
the terminating sum $\displaystyle\app2{-m;b,c}{2b,2c}{x,2-x}$ is identically zero.
\end{lemma}
\begin{proof} 
Here we have a triangular $p\ge0$, $q\ge0$, $p+q\le m$ terminating sum (\ref{appf2}),
as opposed to rectangular $0\le p\le k$, $0\le q\le\ell$ terminating sums frequent in this article.

The $F_2(x,2-x)$ sum is a rational function in $b,c$ of degree at most $2m$. It is enough to show
that this  rational function is zero for infinitely many values of $b,c$. By (\ref{eq:dih12}) and (\ref{eq:dih32}), the $F_2$ sum is zero when $b,c$ are integers $\le -\frac{m+1}2$.
\end{proof}

We have the monodromy group $\ZZ/2\ZZ$ also in the cases $a=-m$ with 
odd  $m>2\max(k,\ell)$ or even $m>2(k+\ell)$. 
Then the hypergeometric functions in (\ref{eq:c2hpg})  themselves terminate, 
while a non-terminating $\hpgo21$ solution is (\ref{eq:dihb}). 
Theorem \ref{th:zerof2} does not apply to the $F_2$ functions 
in (\ref{eq:dih12})--(\ref{eq:dih32}) for odd $m>2\min(k,\ell)$, because generally
\begin{equation}
\lim_{b\to-k} \frac{(-m)_{i+j} (b)_i(c)_j}{i!j!\,(2b)_i(2c)_j}\neq 0 \qquad
\mbox{for}\qquad 2k<i\le m
\end{equation} 
in the triangular sum of Theorem \ref{th:zerof2}, while the same term is taken for zero %$0$
in the rectangular sums in (\ref{eq:dih12})--(\ref{eq:dih32}).

\subsection{Logarithmic solutions}

By Corollary \ref{lm:logc2moodr}, hypergeometric equation $(\ref{hpgde})$ 
with $(\ref{eq:context})$ and $a=-m$ has logarithmic solutions
if and only if exactly one of the conditions (\ref{eq:conda2}) or (\ref{eq:conda3}) holds.
Then one of the equations (\ref{eq:dih12}), (\ref{eq:dih32}) vanishes, 
and the other is a terminating hypergeometric sum. 
The terminating sum is a non-logarithmic solution, obviously. We have
\begin{eqnarray}  \label{eq:nonlogsol}
\frac{\left(\frac{1-m}2\right)_\ell}{\left(\frac12\right)_\ell}\,
\hpg{2}{1}{\!-\frac{m}{2},\ell-\frac{m-1}{2}}{\frac{1}{2}-k}{\,z} = (1-\sqrt{z})^{m}
\app2{\!-m; -k,-\ell}{-2k,-2\ell}{\frac{2\sqrt{z}}{\sqrt{z}-1},\frac2{\sqrt{z}-1}}   \quad
% (1+\sqrt{z})^{m} \app2{\!-m; -k,-\ell}{-2k,-2\ell}{\frac{2\sqrt{z}}{1+\sqrt{z}},\frac2{1+\sqrt{z}}}
\end{eqnarray}
in the case (\ref{eq:conda2}) with $a=-m$, 
since both summands on the right-hand side of (\ref{eq:dih12}) are equal.
We have the same right-hand side expression for the left-hand side of (\ref{eq:dih32}) 
in the case (\ref{eq:conda3}) with $a=-m$.
The $\hpgo21$ function in the vanishing equation  (\ref{eq:dih12}) or (\ref{eq:dih32})
has a logarithmic expression that can be obtained by differentiating 
both sides of the vanishing equation with respect to $a$
and evaluating at $a=-m$.  We have
\begin{eqnarray*}
\frac{d}{da}(1+\sqrt{z})^{-a}=-(1+\sqrt{z})^{-a}\log(1+\sqrt{z})
\end{eqnarray*}
and similarly for $\frac{d}{da}(1-\sqrt{z})^{-a}$. 
On the left-hand side of the vanishing equation (\ref{eq:dih12}) or (\ref{eq:dih32})
we have exactly one (Pochhammer) factor vanishing; to differentiate the product, it is
enough to differentiate only the vanishing Pochhammer symbol.
\begin{lemma} \label{th:dpochh}
The formula for differentiating a Pochhammer symbol is
\begin{equation}
\frac{d}{da}(a)_N=\left\{ \begin{array}{cl}
(a)_N\left(\frac1a+\frac1{a+1}+\ldots+\frac1{a+N-1}\right), & \mbox{if } (a)_N\neq 0,\\
(-1)^a(-a)!(N+a-1)!,& \mbox{if } (a)_N= 0.\end{array}\right.
\end{equation}
\end{lemma}
\begin{proof}
If $(a)_N\neq 0$, we are differentiating a product of $N$ linear functions in $a$.
If $(a)_N=0$ then $a$ is zero or a negative integer. Writing $a=-m+\epsilon$ we have
\begin{equation}
(a)_N=(-m+\epsilon)(1-m+\epsilon)\cdots(-1+\epsilon)\,\epsilon\,(1+\epsilon)\cdots(N-m-1+\epsilon).
\end{equation}
The differentiation is equivalent to dividing out $\epsilon$ and setting $\epsilon=0$.
\end{proof}

Let $\psi(x)$ denote the digamma function, $\psi(x)=\Gamma'(x)/\Gamma(x)$. 
The sum $\frac1a+\frac1{a+1}+\ldots+\frac1{a+N-1}$ can be written as $\psi(a+N)-\psi(a)$ if
$a$ is not zero or a negative integer. 
For an integer $m\ge 0$ let us define the derivative Pochhammer symbol
\begin{equation}
(-m)^\dag_N=\left.\frac{d}{da}(a)_N\right|_{a=-m}.
\end{equation}
By Lemma \ref{th:dpochh} we have
\begin{equation} \label{eq:pochhder}
(-m)^\dag_N=\left\{ \begin{array}{cl}
(-m)_N\big( \psi(m+1-N)-\psi(m+1) \big) & \mbox{if } N\le m,\\
(-1)^mm!\,(N-m-1)!& \mbox{if } N>m.\end{array}\right.
\end{equation}
Note that $(-m)^\dag_N\neq 0$ when $N>0$. On the other hand, $(-m)^\dag_{0}=0$ even if $m=0$. 

When differentiating both sides of (\ref{eq:dih12}) or (\ref{eq:dih32}) with respect to $a$, 
we keep in mind that  $\frac{d}{da}\left(\frac{a}2\right)_{N}$ evaluates to 
$\frac12\left(-\frac{m}2\right)^\dag_{N}$. For even $a=-m$ satisfying (\ref{eq:conda2}), 
we differentiate (\ref{eq:dih32}), use (\ref{eq:nonlogsol}) and  obtain
\begin{eqnarray}  \label{eq:dihlog32}
\frac{\left(\frac{1-m}2\right)_k\!\left(\frac{m}2\right)!\left(k+\ell-\frac{m}2\right)!}
{\left(\frac12\right)_k\left(\frac12\right)_{k+1}\left(\frac12\right)_\ell}\,(-1)^{k+\frac{m}2} z^{k+\frac12}\,
\hpg{2}{1}{k-\frac{m-1}{2},\,k+\ell-\frac{m}{2}+1\,}{\frac{3}{2}+k}{\,z}\hspace{-226pt}\nonumber\\
\equal \frac{\left(\frac{1-m}2\right)_\ell}{\left(\frac12\right)_\ell}\,
\log\frac{1+\sqrt{z}}{1-\sqrt{z}}\;
\hpg{2}{1}{-\frac{m}{2},\,\ell-\frac{m-1}{2}}{\frac{1}{2}-k}{\,z} \nonumber\\
&& +(1-\sqrt{z})^m\,
\appm\dag2{-m; -k,-\ell}{-2k,-2\ell}{\frac{2\sqrt{z}}{\sqrt{z}-1},\frac2{1-\sqrt{z}}}\nonumber\\
&& -(1+\sqrt{z})^m\,
\appm\dag2{-m; -k,-\ell}{-2k,-2\ell}{\frac{2\sqrt{z}}{1+\sqrt{z}},\frac2{1+\sqrt{z}}}. \qquad
\end{eqnarray}
Here the $F^\dag_2$ functions are rectangular sums of $(k+1)(\ell+1)$ terms, defined as 
the $F_2$ sums in (\ref{appf2}) but with each Pochhammer symbol $(-m)_{i+j}$ replaced
by the derivative $(-m)^\dag_{i+j}$ following (\ref{eq:pochhder}). %In particular, 
For odd $a=-m$ satisfying (\ref{eq:conda2}), the differentiation of (\ref{eq:dih32})
 gives the same right-hand side, but the left-hand side is
\begin{equation} \label{eq:dihlog32a}
\hspace{-6pt}
\frac{\left(\frac{m-1}2\right)!\left(k-\frac{m+1}2\right)!\left(-\frac{m}2\right)_{k+\ell+1}}
{\left(\frac12\right)_k\left(\frac12\right)_{k+1}\left(\frac12\right)_\ell}\,(-1)^{k+\frac{m-1}2}z^{k+\frac12}\,
\hpg{2}{1}{\!k-\frac{m-1}{2},k+\ell-\frac{m}{2}+1}{\frac{3}{2}+k}{\,z}.
\end{equation}
For odd $a=-m$ satisfying (\ref{eq:conda3}), % satisfying $2k<m<2\ell$, 
differentiation of (\ref{eq:dih12}) gives 
\begin{eqnarray}   \label{eq:dihlog12}
\frac{\left(\frac{m-1}2\right)!\left(\ell-\frac{m+1}2\right)!}{\left(\frac12\right)_\ell}\,(-1)^{\frac{m-1}2}
\hpg{2}{1}{-\frac{m}{2},\,\ell-\frac{m-1}{2}}{\frac{1}{2}-k}{\,z}=\hspace{-220pt}\nonumber\\
&&\frac{\left(\frac{1-m}2\right)_k\!\left(-\frac{m}2\right)_{k+\ell+1}}
{\left(\frac12\right)_k\left(\frac12\right)_{k+1}\left(\frac12\right)_\ell}\,(-1)^{k}\,z^{k+\frac12}
\log\frac{1+\sqrt{z}}{1-\sqrt{z}}\;
\hpg{2}{1}{\!k-\frac{m-1}{2},k+\ell-\frac{m}{2}+1}{\frac{3}{2}+k}{\,z}\nonumber\\
&&+(1+\sqrt{z})^m\,
\appm\dag2{-m; -k,-\ell}{-2k,-2\ell}{\frac{2\sqrt{z}}{1+\sqrt{z}},\frac2{1+\sqrt{z}}}\nonumber\\
&&+(1-\sqrt{z})^m\,
\appm\dag2{-m; -k,-\ell}{-2k,-2\ell}{\frac{2\sqrt{z}}{\sqrt{z}-1},\frac2{1-\sqrt{z}}}.
\end{eqnarray}
As we will see in Section \ref{sec:trigon}, nice expressions for our dihedral $\hpgo21$ functions
are obtained under the trigonometric substitution $z\mapsto-\tan^2x$.  
In particular, one may recognize
\begin{equation}
\arctan x = \frac1{2i}\,\log\frac{1+ix}{1-ix}.
%\qquad \arcsin x = \frac1{i}\,\log\left(\sqrt{1-x^2}+ix\right).
\end{equation}

\subsection{Degeneration of $24$ Kummer's solutions}
\label{sc:degenerate}

As observed in Remark \ref{kummer24}, the structure of 24 Kummer's solutions
degenerates  if the hypergeometric equation has terminating or logarithmic solutions.
Particurlarly, not all 24 Kummer's solutions are well defined or distinct 
when logarithmic solutions are present.

In the case of the monodromy group $\ZZ/2\ZZ$, the degenerate structure
is described in \cite[Section 7]{degeneratehpg}.
The two solutions in (\ref{eq:c2hpg}) can be represented 
by terminating $\hpgo21$ sums in 6 ways
(with any of the arguments $z,z/(z-1),1/z,1/(1-z),1-1/z,1-z$)
and by non-terminating $\hpgo21$ series (around $z=0$ or $z=\infty$) in 4 ways. 
The remaining four of the 24 Kummer's solutions represent, up to a power factor,
the function in (\ref{eq:diha}). Terminating $F_3$ expression (\ref{eq:diha}) holds, 
but that function %the function in (\ref{eq:diha}) 
can be expressed as a linear combination of two terminating $\hpgo21$ solutions.
For example, the following $z\mapsto1-z$ version of \cite[(43)]{degeneratehpg}
can be used:
\begin{eqnarray} \label{eq:cycl2}
(1-z)^{n+m+1}\hpg21{a+n+1,\,m+1}{n+m+2}{1-z} \!\!\equal\!
\frac{(n+1)_{m+1}}{(-a)_{m+1}}\,
\hpg{2}{1}{-n,\,a-m}{1+a}{\,z}  %\hspace{-12pt} 
\nonumber\\ && \hspace{-29pt} +
\frac{(m+1)_{n+1}}{(a)_{n+1}}
z^{-a}\hpg{2}{1}{-m,-a-n}{1-a}{z}\!,   \qquad
\end{eqnarray}
with
\begin{eqnarray*}
& (n,m,a)\mapsto \left(k-\frac{m+1}2,\ell-\frac{m+1}2,-k-\frac12\right), &\qquad
\mbox{for odd } m<2\min(k,\ell), \\ %\mbox{if (\ref{eq:conda1}) holds},\\ 
& (n,m,a)\mapsto \left(\frac{m-1}2-\ell,\frac{m-1}2-k,-k-\frac12\right), &\qquad
\mbox{for odd } m>2\max(k,\ell),\\ 
& (n,m,a)\mapsto \left(\frac{m}2,\frac{m}2-k-\ell-1,-k-\frac12\right), &\qquad
\mbox{for even } m>2(k+\ell).\;
\end{eqnarray*}

The logarithmic case with terminating solutions is described in \cite[Section 6]{degeneratehpg}.
There are 20 distinct hypergeometric series, or less if $m=k+\ell$. 
Among them, there are 8 terminating and 4 nonterminating hypergeometric series 
representing the non-logarithmic terminating solution, as in (\ref{eq:nonlogsol}). 
The remaining 8 (in general) Kummer's  solutions are non-terminating series
around $z=0$ or $z=\infty$ and represent two different functions, in particular
the logarithmic solution among (\ref{eq:dihlog32})--(\ref{eq:dihlog12}). 
The logarithmic solution can be expressed following \cite[Theorem 6.1]{degeneratehpg} as well.
Formulas \cite[(36) and (38)]{degeneratehpg} give the generic identity
\begin{eqnarray}  \label{eq:logdih32}
\hspace{-12pt}
\frac{(-1)^{m+1}z^{n+m+1-a}}{\left(1-a\right)_{n+m+1}}
\hpg{2}{1}{m+1-a,n+m+1}{n+m+2-a}{\,z}  \hspace{-150pt}  \nonumber  \\
&=& \frac{(m+1-a)_n}{n!\,(n+m)!}\,\hpg{2}{1}{-n,\,a}{a-n-m}{\,z}\left(\log(1-z)+\pi\cot\pi a\right) \nonumber \\
&& -\frac{(z-1)^{-m}}{\left(1-a\right)_{m}}
\sum_{j=0}^{m-1}\frac{\left(a-m\right)_j(m-j-1)!}
{\left(n+m-j\right)!\,j!}(1-z)^j\nonumber\\
&& \!+\! \sum_{j=0}^{n} \!\frac{\left(a\right)_j
\left(\psi\!\left(a\!+\!j\right) + \psi\!\left(n\!-\!j\!+\!1\right)
- \psi(m\!+\!j\!+\!1) - \psi(j\!+\!1)\right)}{(m+j)!\left(n-j\right)!\,j!}(z\!-\!1)^j\nonumber\\
&& +(-1)^{n}\sum_{j=n+1}^{\infty}
\frac{\left(a\right)_j\left(j-n-1\right)!}{(m+j)!\,j!}(1-z)^j,
\end{eqnarray}
where we should substitute
\begin{eqnarray*}
& (n,m,a)\mapsto \left(\frac{m}2,k+\ell-m,\ell-\frac{m-1}2\right), &\qquad
\mbox{if $m$ is even, } 0\le m\le k+\ell, \\ %\mbox{if (\ref{eq:conda1}) holds},\\ 
& (n,m,a)\mapsto \left(k+\ell-\frac{m}2,m-k-\ell,\frac{m+1}2-k\right), &\qquad
\mbox{if $m$ is even, } k+\ell\le m\le2(k+\ell), \\ %\mbox{if (\ref{eq:conda1}) holds},\\ 
& (n,m,a)\mapsto \left(\frac{m-1}2-\ell,k+\ell-m,-\frac{m}2\right), &\qquad
\mbox{if $m$ is odd, } 2\min(k,\ell)<m<2k,\\
& \hspace{-16pt} (n,m,a)\mapsto \left(\frac{m-1}2-k,k+\ell-m,k+\ell-\frac{m}2+1\right), &\qquad
\mbox{if $m$ is odd, } 2\min(k,\ell)<m<2\ell,
\end{eqnarray*}
to get the same $\hpgo21$ function on the left-hand side as in (\ref{eq:dihlog32})--(\ref{eq:dihlog12}).
The $\cot\pi a$ term is zero when $a$ is a half-integer in (\ref{eq:logdih32}). 
These formulas do not contain double sums or square roots, but the last term in  (\ref{eq:logdih32})
is a non-terminating series.

\begin{remark} \rm \label{rm:degenerate}
As observed in Section \ref{sec:easycase}, the terminating 
$F_2$ or $F_3$ sums in (\ref{eq:diha})--(\ref{eq:dih32})
become terminating $\hpgo21$ sums in the case $k=0$ or $\ell=0$. 
If, additionally, $a$ is an integer, the structure of 24 Kummer's solutions 
for the quadratically transformed equation degenerates further.
If $a=-m$ with $0\le m\le 2k$, there are logarithmic solutions by Corollary \ref{lm:logc2moodr}. 
We have to apply then \cite[Section 9]{degeneratehpg} with \mbox{$(n,m,\ell)=(|k-m|,|k-m|,m)$}
to the transformed $\hpgo21$ sums. There are at most 10 different terminating 
$\hpgo21$ solutions, all representing the same elementary solution of the transformed equation.

If $a$ is an integer and $|a+k|>k$, the dihedral group degenerates to $\ZZ/2\ZZ$, while the monodromy group of the quadratically transformed equation is trivial.  By \cite[Section 8]{degeneratehpg},
the transformed equation has 3 solutions like in (\ref{eq:termlz}), %--(\ref{eq:termlzz}), 
each representable by 6 terminating and 2 non-terminating $\hpgo21$ series.
\end{remark}

\section{Trigonometric expressions}
\label{sec:trigon}

Formulas (\ref{er:triga})--(\ref{er:trigb}) show attractive trigonometric expressions for the simplest
dihedral functions. Analogous trigonometric expressions are possible for all dihedral functions,
if only due to contiguous relations. The literature appears to give only trigonometric expressions
for the simplest dihedral functions. In particular \cite[15.1]{abrostegun}, \cite[2.8]{bateman},
here are fractional-linear transformations of  (\ref{er:triga})--(\ref{er:trigb}):
\begin{eqnarray} 
\hpg21{\,\frac{1+a}2,\,\frac{1-a}2\,}{\frac12}{\sin^2 x} \equal \frac{\cos ax}{\cos x},\\ 
\hpg21{\frac{2+a}2,\frac{2-a}2}{\frac32}{\sin^2 x} \equal \frac{2\sin ax}{a\sin 2x},\\
\hpg21{\,\frac{a}2,\,\frac{a+1}2\,}{\frac12}{-\tan^2 x} \equal \cos ax\,\cos^a x,\\ 
\hpg21{\frac{a+1}2,\frac{a+2}2}{\frac32}{-\tan^2 x} \equal \frac{\sin ax}{a\sin x}\,\cos^{a+1}x.
\end{eqnarray}
According to the latter two formulas, we should use the argument substitution $z\mapsto-\tan^2x$ 
in our main expressions. Nicer expressions are obtained after fractional-linear transformation
(\ref{flinear3}), hence with the argument $\sin^2 x$.

\subsection{General trigonometric expressions}

The following theorem gives trigonometric versions of formulas  (\ref{eq:dih12})--(\ref{eq:dih32})
after Pfaff's transformation (\ref{flinear3}).
Subsequently, we discuss here simplification of trigonometric formulas. 
Trigonometric modification of formula (\ref{eq:diha}) is discussed
in Subsection \ref{sc:hyperbolic}.
%, and trigonometric forms
%of degenerate or invariant expressions for Sections \ref{sec:symsq}, \ref{adihedral}
% and \ref{sec:degenerate}. 
\begin{theorem}
Let us denote
\begin{equation}
\Upsilon_{p,q}^{a,k,\ell}(x)
:=\frac{2^{p+q}\,(-k)_p\,(-\ell)_q\,(a)_{p+q}}{(-2k)_p\,(-2\ell)_q\;p!\,q!}\,\sin^px\,\cos^q x.
\end{equation}
Then
\begin{eqnarray}  \label{eq:diht12}
\frac{\left(\frac{a+1}2\right)_\ell}{\left(\frac12\right)_\ell}\;
\hpg{2}{1}{\,\frac{a}{2},\,-\frac{a}{2}-k-\ell\,}{\frac{1}{2}-k}{\sin^2 x} \!\!\!\equal\!\!
\sum_{p=0}^{k} \sum_{q=0}^{\ell}
%\frac{2^{p+q}\,(-k)_p\,(-\ell)_q\,(a)_{p+q}}{(-2k)_p\,(-2\ell)_q\;p!\,q!}\,\sin^px\,\cos^q x\,
\Upsilon_{p,q}^{a,k,\ell}(x)\, {\textstyle \cos\left(ax+(p+q)x-\frac{\pi}2p\right), }
\qquad\\ \label{eq:diht32} %\hspace{-11pt}
\frac{\left(\frac{a+1}2\right)_k\left(\frac{a}2\right)_{k+\ell+1}\,\sin^{2k+1} x}
{\left(\frac12\right)_k\left(\frac12\right)_{k+1}\left(\frac12\right)_\ell}\,
\hpg{2}{1}{\!\frac{a+1}{2}+k,\frac{1-a}{2}-\ell}{\frac{3}{2}+k}{\sin^2x} 
\hspace{-99pt} & \nonumber\\  
\!\!\!\equal\!\! \sum_{p=0}^{k} \sum_{q=0}^{\ell}
\Upsilon_{p,q}^{a,k,\ell}(x)\, {\textstyle \sin\left(ax+(p+q)x-\frac{\pi}2p\right). }
\end{eqnarray}
\end{theorem}
\begin{proof} After %application of 
fractional-linear transformation (\ref{flinear3})
and the substitution \mbox{$z\mapsto-\tan^2x$} in (\ref{eq:dih12}) we get
\begin{eqnarray*} \hspace{-6pt}
\frac{\left(\frac{a+1}2\right)_\ell}{\left(\frac12\right)_\ell}\,
\hpg{2}{1}{\frac{a}{2},-\frac{a}{2}-k-\ell}{\frac{1}{2}-k}{\sin^2x}\equal\frac{\exp(-iax)}2\,
\app2{a; -k,-\ell}{-2k,-2\ell}{\frac{2i\sin x}{\exp ix},\frac{2\cos x}{\exp ix}}\nonumber\\
&& +\frac{\exp iax}2\,
\app2{a; -k,-\ell}{-2k,-2\ell}{\frac{-2i\sin x}{\exp(-ix)},\frac{2\cos x}{\exp(-ix)}}.
\end{eqnarray*}
%\vspace{1pt}
We use Euler's formula $\exp ix=\cos x+i\sin x$ as well. We add the two double sums term by term,
and after the identification $\pm2i=2\exp(\pm i\frac{\pi}2)$ we get the double sum in (\ref{eq:diht12}).
Note that $\Upsilon_{p,q}^{a,k,\ell}(x)$
is the generic summand of $\displaystyle\app2{a;\,-k,-\ell}{-2k,-2\ell}{2\sin x,2\cos x}$.

Formula (\ref{eq:diht32}) follows similarly. Withal, the substitution \mbox{$z\mapsto-\tan^2x$}
into $(-1)^kz^{k+1/2}$ gives $i\sin^{2k+1}x/\cos^{2k+1}x$.
\end{proof}

The right-hand side of (\ref{eq:diht12}) can be written as
\begin{equation} \label{eql:dih12}
P_{k,\ell}(x)\,\cos ax+Q_{k,\ell}(x)\,\sin ax,
\end{equation}
with 
\begin{eqnarray*}
P_{k,\ell}(x) &=& \sum_{p=0}^{k} \sum_{q=0}^{\ell}
\Upsilon_{p,q}^{a,k,\ell}(x)\, {\textstyle \cos\left(\frac{\pi}2p-(p+q)x\right), }\\
Q_{k,\ell}(x) &=& \sum_{p=0}^{k} \sum_{q=0}^{\ell}
\Upsilon_{p,q}^{a,k,\ell}(x)\, {\textstyle \sin\left(\frac{\pi}2p-(p+q)x\right). }
\end{eqnarray*}
The right-hand side of (\ref{eq:diht32}) can be then written as
\begin{equation}  \label{eql:dih32}
P_{k,\ell}(x)\,\sin ax-Q_{k,\ell}(x)\,\cos ax,
\end{equation}
and we have
\begin{equation} \label{eq:intkl} \textstyle
P_{\ell,k}(x)=P_{k,\ell}(\frac{\pi}2-x), \qquad Q_{\ell,k}(x)=-Q_{k,\ell}(\frac{\pi}2-x).
\end{equation}
Generally, here are the first few terms of the double sums $P_{k,\ell}(x)$, $Q_{k,\ell}(x)$:
\begin{eqnarray*}
P_{k,\ell}(x) \!\equal\! 1+a\cos^2x+\frac{a(a+1)(\ell-1)}{2\ell-1}\cos^2x\cos2x+
\frac{a(a+1)(a+2)(\ell-2)}{3(2\ell-1)}\cos^3x\cos3x+\ldots\\
&&+a\sin^2x+\frac{a(a+1)}{2}\sin^22x+
\frac{a(a+1)(a+2)(\ell-1)}{2\ell-1}\sin x\cos^2x\sin3x+\ldots\\
&&-\frac{a(a+1)(k-1)}{2k-1}\sin^2x\cos2x-
\frac{a(a+1)(a+2)(k-1)}{2k-1}\sin^2 x\cos x\cos3x+\ldots\\
&&-\frac{a(a+1)(a+2)(k-2)}{3(2k-1)}\sin^3x\sin3x-\ldots\\
%\frac{a(a+1)(a+2)(a+3)(k-2)}{3(2k-1)}\sin^3x\cos x\sin4x+
&&+\ldots,\\
Q_{k,\ell}(x) \!\equal\! 0-\frac{a}2\sin2x-\frac{a(a+1)(\ell-1)}{2\ell-1}\cos^2x\sin2x-
\frac{a(a+1)(a+2)(\ell-2)}{3(2\ell-1)}\cos^3x\sin3x-\ldots\\
&&+\frac{a}2\sin2x+\frac{a(a+1)}{4}\sin4x
+\frac{a(a+1)(a+2)(\ell-1)}{2\ell-1}\sin x\cos^2x\cos3x+\ldots\\
&&+\frac{a(a+1)(k-1)}{2k-1}\sin^2x\sin2x+
\frac{a(a+1)(a+2)(k-1)}{2k-1}\sin^2 x\cos x\sin3x+\ldots\\
&&-\frac{a(a+1)(a+2)(k-2)}{3(2k-1)}\sin^3x\cos3x-\ldots\\
&&-\ldots.
\end{eqnarray*}
Evidently, significant simplification in the double sums $P_{k,\ell}(x)$, $Q_{k,\ell}(x)$ is possible
if we collect terms to the same Pochhammer products $(a)_{p+q}$. If both $k\ge 3$ and $\ell\ge 3$,
we have
\begin{eqnarray*}
P_{k,\ell}(x) \!\equal\! 1+a+\frac{a(a+1)}2\,\big(1-B_{k,\ell}(x)\cos2x\big)+
\frac{a(a+1)(a+2)}{6}\,\big(1-3B_{k,\ell}(x)\cos2x\big)+\ldots,\\
Q_{k,\ell}(x) \!\equal\! \frac{a(a+1)}2\,B_{k,\ell}(x)\sin2x+
\frac{a(a+1)(a+2)}{2}\,B_{k,\ell}(x)\sin2x+\ldots,
\end{eqnarray*}
where 
$$B_{k,\ell}(x)=\frac{\cos^2x}{2\ell-1}-\frac{\sin^2x}{2k-1}.$$ %\vspace{2pt}

Simplification in the symmetric case $k=\ell$ is particularly significant,
as could be expected from relation (\ref{eq:symmkl}) to the case $\ell=0$. 
Quadratic transformation (\ref{quadr1}) can be written as 
\begin{equation} \label{eq:quadrsin}
\hpg{2}{1}{a,\,b}{\frac{a+b+1}{2}}{\sin^2 x} =
\hpg{2}{1}{\frac{a}{2},\,\frac{b}{2}}{\frac{a+b+1}{2}}{\sin^2 2x}.
\end{equation}
This relates the case $\ell=k$ with $a\mapsto 2a$ of (\ref{eq:diht12}) 
to the case $\ell=0$ with $x\to 2x$ of (\ref{eq:diht12}). To use (\ref{eq:quadrsin})
similarly with (\ref{eq:diht32}), one has to apply Euler's transformation (\ref{flinear1}) to 
the $\hpgo21$ function in (\ref{eq:diht32}) first.

For the dihedral functions with small $k$, $\ell$
nice expressions are obtained if we make the substitution
$a\mapsto \cc-k-\ell$ and write the right-hand sides of (\ref{eq:diht12}) and (\ref{eq:diht32})
as, respectively,
\begin{eqnarray} \label{eq:}
{\bf P}_{k,\ell}(x)\,\cos \cc x+{\bf Q}_{k,\ell}(x)\,\sin \cc x,\\
{\bf P}_{k,\ell}(x)\,\sin \cc x-{\bf Q}_{k,\ell}(x)\,\cos \cc x.
\end{eqnarray}
Explicitly, we have the rotation transformation
\begin{equation} \label{eq:rottr}
\begin{array}{rcl}
{\bf P}_{k,\ell}(x) &=& P_{k,\ell}(x)\cos(k+\ell)x-Q_{k,\ell}\sin(k+\ell)x,\\
{\bf Q}_{k,\ell}(x) &=& P_{k,\ell}(x)\sin(k+\ell)x+Q_{k,\ell}\cos(k+\ell)x.
\end{array}
\end{equation}
A bit tedious translation of (\ref{eq:intkl}) shows the following effect of
interchange of the indices $k,\ell$:
\begin{eqnarray} \label{eq:intkla}
{\bf P}_{\ell,k}(x) &=& \left\{ \begin{array}{cl}
(-1)^j\,{\bf P}_{k,\ell}(\frac{\pi}2-x), & \mbox{if $k+\ell=2j$ even},\\
(-1)^j\,{\bf Q}_{k,\ell}(\frac{\pi}2-x), & \mbox{if $k+\ell=2j+1$ odd},
\end{array} \right.\\    \label{eq:intklb}
{\bf Q}_{\ell,k}(x) &=& \left\{ \begin{array}{cl}
(-1)^{j-1}\,{\bf Q}_{k,\ell}(\frac{\pi}2-x), & \mbox{if $k+\ell=2j$ even},\\
(-1)^j\,{\bf P}_{k,\ell}(\frac{\pi}2-x), & \mbox{if $k+\ell=2j+1$ odd}.
\end{array} \right.
\end{eqnarray}
But here are attractive evaluation formulas after the substitution $a\mapsto \cc-k-\ell$:
\begin{eqnarray}
\hpg{2}{1}{\frac{-1+\cc}{2},\,\frac{-1-\cc}{2}}{\frac{1}{2}}{\sin^2 x} &=&
\cos x \cos\cc x+\sin x\,\frac{\sin\cc x}{\cc},\\
\hpg{2}{1}{\frac{-1+\cc}{2},\,\frac{-1-\cc}{2}}{-\frac{1}{2}}{\sin^2 x} &=&
\cos x \cos\cc x+\cc\,\sin x\,\sin\cc x,\\
\hpg{2}{1}{\,\frac{\cc}{2},\,-\frac{\cc}{2}}{\frac{3}{2}}{\sin^2 x} &=&
\frac{\cc\,\cos x \sin\cc x-\sin x\,\cos\cc x}{(\cc^2-1)\,\sin x},\\
\hpg{2}{1}{\frac{2+\cc}{2},\,\frac{2-\cc}{2}}{\frac{5}{2}}{\sin^2 x} &=&
3\,\frac{\cos x \sin\cc x-\cc\,\sin x\,\cos\cc x}
{\cc\,(\cc-1)\,(\cc+1)\,\sin^3 x},\\
\hpg{2}{1}{\frac{-2+\cc}{2},\,\frac{-2-\cc}{2}}{-\frac{1}{2}}{\sin^2 x} &=&
\cos 2x \cos\cc x+\frac{\cc}2\,\sin 2x\,\sin\cc x,\\
\hpg{2}{1}{\frac{3+\cc}{2},\,\frac{3-\cc}{2}}{\frac{5}{2}}{\sin^2 x} &=&
3\,\frac{\cos 2x \sin\cc x-\frac12\,\cc\,\sin 2x\,\cos\cc x}
{\cc\,(\cc-2)\,(\cc+2)\,\sin^3 x}.
\end{eqnarray}
Some further expressions:
\begin{eqnarray*}
{\bf P}_{0,2}(x)&=&\frac{\cc^2\cos^2x-2\sin^2x-1}3, 
\qquad {\bf Q}_{0,2}(x)=\frac{\cc\,\sin 2x}2,\\
{\bf P}_{1,2}(x)&=&\frac{\cc(\cc-2)}3\,\cos x\,(2\cos2x-1), \\
{\bf Q}_{1,2}(x)&=&\frac{\cc-2}3\,\sin x\,(\cc^2\cos^2x-3\sin^2x),\\
{\bf P}_{2,2}(x)&=&\frac{(\cc-1)(\cc-3)}9\left(\cos 4x+2-\frac{\cc^2\sin^22x}4\right), \\
{\bf Q}_{2,2}(x)&=&\frac{\cc\,(\cc-1)\,(\cc-3)}{12}\,\sin 4x,\\
{\bf P}_{0,3}(x)&=&\frac{\cc\cos x\left(\cc^2\cos^2 x-11\sin^2 x-4\right)}{15}, \\
{\bf Q}_{0,3}(x)&=&\frac{\sin x\,(2\cc^2\sin^2 x-2\sin^2 x-3)}{5},\\
{\bf P}_{1,3}(x)&=&\frac{\cc-3}{15}\,\left( \cc^2\cos^2x\,(1-7\sin^2 x)
+\cos 4x+16\sin^2 x-2\right),\\
{\bf Q}_{1,3}(x)&=&\frac{\cc(\cc-3)}{30}\,\sin 2x\,(\cc^2\cos^2x-14\sin^2x-1), \\
{\bf P}_{3,3}(x)&=&\frac{(\cc-1)(\cc-3)(\cc-5)}{75}\,\cos 2x
\left(\cos 4x+4-\frac{\cc^2\sin^22x}2\right), \\
{\bf Q}_{3,3}(x)&=&\frac{\cc\,(\cc-1)\,(\cc-3)\,(\cc-5)}{900}\,\sin 2x
\left(11\cos 4x+19-\frac{\cc^2\sin^22x}2\right),\\
{\bf P}_{0,4}(x)&=&\frac{\cc^4\cos^4x}{105}
-\frac{\cc^2\cos^2 x\,(7\sin^2x+2)}{21}+\frac{\cos4x+32\sin^2x+2}{35}, \\
{\bf Q}_{0,4}(x)&=&\frac{\cc\,\sin 2x\,(2\cc^2\cos^2 x-10\sin^2 x-11)}{42}.
 \end{eqnarray*}
 
 \subsection{Apparent Pad\'e approximation properties}
 
 Because of the $O(x^{2k+1})$ solution in (\ref{eq:diht32}), the introduced trigonometric polynomials
 have properties resembling simultaneous Pad\'e approximation. We have
 \begin{equation}
 \frac{Q_{k,\ell}(x)}{P_{k,\ell}(x)}=\tan a x+O\left(x^{2k+1}\right), \qquad
  \frac{{\bf Q}_{k,\ell}(x)}{{\bf P}_{k,\ell}(x)}=\tan \cc x+O\left(x^{2k+1}\right).
 \end{equation}
and by the symmetry (\ref{eq:intkl}) or  (\ref{eq:intkla})--(\ref{eq:intklb})  of $k,\ell$,
 \begin{equation}
 \frac{Q_{k,\ell}(\frac{\pi}2-x)}{P_{k,\ell}(\frac{\pi}2-x)}=-\tan a x+O\left(x^{2\ell+1}\right), \quad
\left(\frac{{\bf P}_{k,\ell}(\frac{\pi}2-x)}{{\bf Q}_{k,\ell}(\frac{\pi}2-x)} \right)^{\!\varepsilon}
=\varepsilon\tan \cc x+O\left(x^{2\ell+1}\right).
 \end{equation}
Here $\varepsilon=(-1)^{k+\ell-1}$. It is however not clear what characterizes
the trigonometric polynomials uniquely in the approximation context.
From the degenerate cases with vanishing Pochhammer factors 
on the right-hand sides of  (\ref{eq:diht12}) or (\ref{eq:diht32}) 
we get the following vanishing  or exact approximation properties:
%with $K_1=\min(k,\ell)$, $K_2=\ell-k$, $L=k+\ell$:
\begin{eqnarray} \label{eq:pqdegen}
P_{k,\ell}(x)=0,\quad Q_{k,\ell}(x)=0, &\qquad& \mbox{if (\ref{eq:conda1}) holds}, \nonumber\\
\frac{Q_{k,\ell}(x)}{P_{k,\ell}(x)} =\tan ax, \quad && \mbox{if (\ref{eq:conda2}) holds}, \nonumber\\
\frac{P_{k,\ell}(x)}{Q_{k,\ell}(x)} =-\tan ax, && \mbox{if (\ref{eq:conda3}) holds},\\
{\bf P}_{k,\ell}(x)=0,\quad {\bf Q}_{k,\ell}(x)=0, &&
\mbox{if (\ref{eq:conde1}) holds}, \nonumber\\
\frac{{\bf Q}_{k,\ell}(x)}{{\bf P}_{k,\ell}(x)} =\tan \cc x, \quad && 
\mbox{if (\ref{eq:conde2}) holds}, \nonumber\\
\frac{{\bf P}_{k,\ell}(x)}{{\bf Q}_{k,\ell}(x)} =-\tan \cc x, && 
\mbox{if (\ref{eq:conde3}) holds},\nonumber
\end{eqnarray}
and similarly with $P_{k,\ell}(\frac{\pi}2-x),Q_{k,\ell}(\frac{\pi}2-x)$ or 
${\bf P}_{k,\ell}(\frac{\pi}2-x),{\bf Q}_{k,\ell}(\frac{\pi}2-x)$.
These conditions do not determine the trigonometric polynomials uniquely either,
because ${\bf P}_{k,\ell}(x),{\bf Q}_{k,\ell}(x)$ are not
the shifted $a\mapsto\cc-k-\ell$ versions of $P_{k,\ell}(x),Q_{k,\ell}(x)$.

\subsection{Contiguous relations}

Contiguous relations for the considered $\hpgo21$ functions translate into
the following recurrence relations:
\begin{eqnarray}
P_{k,\ell}(x) \!\!\equal\!\! \frac{a+2k-1}{2k-1}\sin^2x\,P_{k-1,\ell}(x)
+\frac{a+2\ell-1}{2\ell-1}\cos^2x\,P_{k,\ell-1}(x), \qquad\quad \\
\!\!\left(\frac{a}2+k+\ell+1\right)\!P_{k,\ell}(x) \!\!\equal\!\! %\textstyle
\left(\,k+\frac12\,\right)P_{k+1,\ell}(x)+\left(\,\ell+\frac12\,\right)P_{k,\ell+1}(x),\\
\left(\,\ell+\frac12\,\right)P_{k,\ell+1}(x) \!\!\equal\!\!
\left(\,\ell+\frac12+(a+k+\ell)\cos^2x\,\right)P_{k,\ell}(x) \nonumber\\
&&-\frac{a+2\ell-1}{2\ell-1}\left(\,\frac{a}2+k+\ell\,\right)\cos^2x\,P_{k,\ell-1}(x).
\end{eqnarray}
The recurrence relations for $Q_{k,\ell}$'s are the same.
Three-term recurrences %relations
for the trigonometric polynomials  ${\bf P}_{k,\ell}(x)$, ${\bf Q}_{k,\ell}(x)$ are less clear, 
because some shifts in $k,\ell$ lead to half-integer shifts in the upper parameters
of the $\hpgo21$ functions, and transformation (\ref{eq:rottr}) mixes ${\bf P}_{k,\ell}$'s
with ${\bf Q}_{k,\ell}$'s. Since we formally have 
\mbox{${\bf Q}_{0,0}(x)=0$}, \mbox{${\bf Q}_{0,-1}(x)=0$}, 
three-term recurrences between these polynomials with minimal shifts in $k,\ell$ 
should not be expected. Here are some consequences of the contiguous relations, 
nevertheless:
\begin{eqnarray}
\left(2k\cos^2x+2\ell\sin^2x+1\right) {\bf P}_{k,\ell}(x) \!\equal\!
\frac{2k+1}{2\ell-1}\left(\cc-k+\ell-1\right)\cos^2x\,{\bf P}_{k+1,\ell-1}(x) \nonumber \\
&&+\frac{2\ell+1}{2k-1}\left(\cc+k-\ell-1\right)\sin^2x\,{\bf P}_{k-1,\ell+1}(x), \qquad\\
%\frac{\beta-k+\ell-1}{(2\ell-1)(2\ell+1)}\cos^2x\,{\bf P}_{k+1,\ell-1}(x) &= &
\frac{\sin^2x}{2k+1}\,{\bf P}_{k,\ell}(x)+\frac{{\bf P}_{k+1,\ell+1}(x)}{\cc-k-\ell-1}
\!\equal\! \frac{\cc-k+\ell-1}{(2\ell-1)(2\ell+1)}\cos^2x\,{\bf P}_{k+1,\ell-1}(x), \\
(2\ell-1)\,{\bf P}_{k,\ell}(x) \!\equal\! %\textstyle
\frac{\cc^2-(k+\ell)^2}{(2k-1)(2k+1)}(\cc-k-\ell+1)\sin^2x\,{\bf P}_{k-1,\ell-1}(x)
\nonumber \\
&& +(\cc-k+\ell-1)\,{\bf P}_{k+1,\ell-1}(x).
\end{eqnarray}

\subsection{A hyperbolic trigonometric version of $(\ref{eq:diha})$}
\label{sc:hyperbolic}

A trigonometric version of formula (\ref{eq:diha}) is pointless,
because the substitution \mbox{$z\mapsto-\tan^2x$} does not give a meaningful 
$\hpgo21$ argument for a local expansion. We may however use 
the hyperbolic tangent substitution $z\mapsto \tanh^2x$
(plus $a\mapsto \cc-k-\ell$ for a little neatness), 
and consider the approximation $x\to+\infty$:
\begin{eqnarray} 
\label{eq:dihta} %\hspace{-9pt}
\hpg{2}{1}{\frac{\cc-k-\ell}{2},\,\frac{\cc-k+\ell+1}{2}}{\cc+1}{\frac{1}{\cosh^2 x}} \!\!\equal\!
2^\cc\cosh^{\cc} x\,\left(\cosh \cc x-\sinh \cc x \right)\,\tanh^k x
\times\nonumber\\ 
&& \hspace{-30pt} 
\app3{k+1,\ell+1; -k,-\ell}{\cc+1}{\frac{1-\coth x}{2},\frac{1-\tanh x}{2}}\!. \qquad
\end{eqnarray}

Trigonometric expressions for the dihedral functions with $k\ell=0$ or $k=\ell$
could be obtained in relation to the Legendre functions (as described in Remark \ref{rm:legendre}).
Infinite trigonometric series for general Legendre functions are given
in \cite[8.7.1--8.7.2]{abrostegun}.

\subsection{Trigonometric expressions for logarithmic solutions}

The use of trigonometric arguments is instructive for
degenerate or logarithmic dihedral $\hpgo21$ functions as well,
as suggested by formulas (\ref{dihedr4}), (\ref{dihedr5}). 
In the logarithmic cases (\ref{eq:conda2}), (\ref{eq:conda3})
the non-logarithmic solution like (\ref{eq:nonlogsol}) 
in (\ref{eq:dih12}) or (\ref{eq:dih32})  is, respectively,
\begin{equation}
\frac{P_{k,\ell}(x)}{\cos ax}, \qquad\qquad  -\frac{Q_{k,\ell}(x)}{\cos ax},
\end{equation}
due to (\ref{eq:pqdegen}). The logarithmic solution is obtained by differentiating
(\ref{eql:dih32}) or (\ref{eql:dih12})  with respect to $a$. 
Under condition (\ref{eq:conda2}), the left hand side of (\ref{eq:dihlog32}) if $a$ is even,
or (\ref{eq:dihlog32a}) if $a$ is odd, is equal to
\begin{equation}
\frac{x\,P_{k,\ell}(x)}{\cos ax}+\frac{dP_{k,\ell}(x)}{da}\sin ax-\frac{dQ_{k,\ell}(x)}{da}\cos ax,
\qquad\mbox{with }a=-m.
\end{equation} 
Under condition (\ref{eq:conda3}), the left hand side of (\ref{eq:dihlog12}) is equal to
\begin{equation}
\frac{x\,Q_{k,\ell}(x)}{\cos ax}+\frac{dP_{k,\ell}(x)}{da}\cos ax+\frac{dQ_{k,\ell}(x)}{da}\sin ax,
\qquad\mbox{with }a=-m.
\end{equation} 
Similar expressions can be obtained with the more compact trigonometric polynomials
${\bf P}_{k,\ell}(x), {\bf Q}_{k,\ell}(x)$, which depend on the parameter $\cc=a+k+\ell$.

\bibliographystyle{plain}
\bibliography{../../hypergeometric}

\end{document}